\documentclass[16pt]{article}
\usepackage{amsmath,amsfonts,amssymb}
\usepackage{verbatim}
\usepackage{latexsym}
\newtheorem{theorem}{Theorem}[section]
\newtheorem{lemma}[theorem]{Lemma}
\newtheorem{proposition}[theorem]{Proposition}
\newtheorem{definition}[theorem]{Definition}
\newtheorem{corollary}[theorem]{Corollary}

\newcommand\RR{{\Bbb R}}
\newcommand\CC{{\Bbb C}}

\newcommand\ZZ{{\Bbb Z}}

\newcommand\HH{{\Bbb H}}
\newcommand\BB{\mathcal{B}}

\begin{document}
\title{
Continuous Wavelets and Frames\\ on Stratified Lie Groups I
}
\author{Daryl Geller\\
\footnotesize\texttt{{daryl@math.sunysb.edu}}\\
 \thanks{Research supported by the German Academic  Exchange Service (DAAD).
 }
Azita Mayeli \\
\footnotesize\texttt{{mayeli@ma.tum.de}}
}

\maketitle

\begin{abstract}

Let $G$ be a stratified Lie group and  $L$ be the sub-Laplacian on $G$. Let $0 \neq f\in
\mathcal{S}(\RR^+)$.
We show that $Lf(L)\delta$, the distribution kernel of the operator $Lf(L)$, is an admissible
function on $G$.  We also show that, if $\xi f(\xi)$ satisfies Daubechies' criterion, then
$L f(L)\delta$ generates a frame for any sufficiently fine lattice subgroup of $G$.

\end{abstract}

\section{\large{Introduction}}
\label{sec:Introduction }

Let $L$ denote the sub-Laplacian on a stratified group $G$ \cite{FollandStein82}( for instance ,
the Heisenberg group
$\HH^n$).  
If
$\phi\in \mathcal{S}({G})$ and $\int \phi =0$, we say $\phi$ is \textit{admissible } if for some
$c\not=0$,
Calder\'on's reproducing formula:

\begin{align}\label{eq:coroner}
\int_0^\infty \tilde{\phi}_a \ast \phi_a \; a^{-1}da = c\delta.
\end{align}

holds in the sense of tempered distributions,
where $\phi_a(x)=a^{-Q}\phi(a^{-1}x)$, $Q$ is the homogeneous dimension of $G$,
$\widetilde{\phi\textsc{}}(x)= \overline{\phi}(x^{-1})$ and $\delta$ denotes the point mass at
$0\in G$.
(In section 5, we shall show that this definition of ``admissible'' is equivalent to the one
generally
used in wavelet theory.) In section 4, we shall show:

\begin{theorem}\label{T:2}
Let $f$ be a nonzero element of $\mathcal{S}(\RR^+)$. Then $Lf(L)\delta \in \mathcal{S}(G)$ is
admissible.
\end{theorem}

For example, $Le^{-L/2}\delta$ is admissible. (Here $f(L)\delta$ is the distribution kernel of
$f(L)$, so that if $F$ is a Schwartz function, $f(L)F = F*[f(L)\delta]$.)  Up to a constant,
$Le^{-L/2}\delta$ is a very natural generalization of the 
Mexican Hat Wavelet to $G$.  In case $G = \HH^n$, Theorem \ref{T:2} was shown for this
function in Mayeli \cite{Mayelithesis05}.\\

As a corollary of Theorem \ref{T:2}, we shall show in sections 4 and 5:

\begin{corollary}\label{C:1}
(a) There exist admissible $\phi \in \mathcal{S}(G)$ with all moments vanishing.\\
(b) There exist admissible $\phi \in C_c^{\infty}(G)$ with arbitrarily many moments vanishing.
\end{corollary}

In Corollary \ref{C:1} (a) and (b), we will in fact show that $\phi$ can be chosen to have the
form $\phi=Lf(L)\delta$ for some $f \in \mathcal{S}(\RR^+)$. As we will explain at the end
of section 4, Corollary \ref{C:1} improves on Lemmas 1.61 and 1.62 of Folland-Stein
\cite{FollandStein82} for stratified groups. \\

Moreover, we shall show in section 7:  
\begin{theorem}\label{T:3}
Let $\Gamma$ be a lattice subgroup of $G$, and let $f$ again be a nonzero
element of $\mathcal{S} (\RR^+)$.\\ 
(a) If $\xi f(\xi)$ satisfies ``Daubechies' criterion'' then for sufficiently small
$b > 0$, the admissible function
$Lf(L)\delta$ generates a wavelet frame for the lattice $b\Gamma$.  \\
(Note: Daubechies' criterion holds here in particular if 
$f(\xi)$ does not vanish for any $\xi > 0$, or alternatively
if the dilation parameter $a$ is sufficiently close to $1$.)\\
(b) As $a \rightarrow 1$, the ratio of the optimal frame bounds in (a) is 
$1 + O(|a-1|^2 \log|a-1|)$, for sufficiently small $b > 0$.  (Here 
$a$ is again the dilation parameter.)
\end{theorem}

Theorem \ref{T:3} (b) says, in essence, that if $a$ is close to $1$, then the 
frame is ``nearly tight'', and that the
convergence of the ratio of the optimal frame bounds to $1$ is {\em nearly quadratic}
in $|a-1|$.  (Again, $b$ must be sufficiently small, and is chosen after $a$ is chosen.)\\ 

In particular, we shall show that, if one
uses the dilation parameter $a = 2^{1/3}$,
then for all sufficiently small $b > 0$, the  
admissible function $Le^{-L/2}\delta$ generates a 
wavelet frame for $b\Gamma$ which is ``nearly tight'': there are frame bounds
$B_b, A_b$ with $B_b/A_b = 1.0000$ to four significant digits. 
This example shows that $a$ need not be all that close to $1$ for
a nearly tight frame to be obtained in Theorem \ref{T:3} (b); 
$a=2^{1/3}$ is already very good. \\
 
Instead of using $Le^{-L/2}\delta$ as the
admissible function, one could choose $f \in \mathcal{S}(\RR^+)$ so that
$\phi = Lf(L)\delta$ is as in Corollary \ref{C:1}
(a) or (b).  One then obtains nearly tight frames of Schwartz functions with all moments 
vanishing, or nearly tight frames of $C_c^{\infty}$ functions with arbitrarily many
moments vanishing (for suitable $a$ and $b$).\\  

To clarify our terminology in Theorem \ref{T:3}:
\begin{itemize}
\item
$b\Gamma = \{b\gamma: \gamma \in \Gamma\}$; here $b\gamma$, a dilate of 
$\gamma$, is defined in (\ref{eq:dil}) below.
\item
For a fixed dilation parameter $a > 0$,
if $\phi$ is a function on $G$, $j \in \ZZ$ and $\gamma \in \Gamma$, we set
\[\phi_{j,b\gamma}(x) =
a^{-jQ/2} \phi([b\gamma]^{-1}[a^{-j}x]). \]
\item
To say that an $L^2$ function $\psi$ generates a wavelet frame for the lattice
$b\Gamma$ is to say that \newline
$\{\phi_{j,b\gamma}(x) : j \in \ZZ, \gamma \in \Gamma\}$ is a frame.  
\item
To say that a function $g \in \mathcal{S}(\RR^+)$ satisfies Daubechies' criterion
is to say that
\begin{equation}
\label{daubcrit}
A = \inf_{\lambda > 0} \sum_{j=-\infty}^{\infty} |g(a^{2j}\lambda)|^2 > 0. 
\end{equation}

In \cite{Daubechies92}, page 68, Daubechies observes that if $G = \RR$ and 
$\Gamma = \ZZ$, then this is a necessary condtion in Theorem \ref{T:3} (a).  
Here we have put $g(\lambda) = \lambda f(\lambda)$, for $f \in \mathcal{S}(\RR^+)$.  
Then it is easily
seen that the series in (\ref{daubcrit}) converges uniformly on compact subsets of
$(0,\infty)$.  Let $u(\lambda)$ denote the sum of that series; then clearly
$u(a^2 \lambda) = u(\lambda)$ for all $\lambda > 0$.
Consequently, $A$ is the just the minimum of the series for 
$\lambda \in [1,a^2]$.  Thus Daubechies' criterion is equivalent to the nonexistence
of a $\lambda_0 > 0$ such that $g(a^{2j} \lambda_0) = 0$ for all integers $j$.

\end{itemize}

In fact, in Theorem \ref{T:3}, one does not even
need the full force of the assumption that $\Gamma$ is a lattice subgroup; all that
one needs is that $\Gamma$ is a discrete subset of $G$, and that
there is a bounded measurable set $\mathcal{R}$, of positive
measure, such that every $g \in G$ may be written uniquely in the form 
$g = x\gamma$ with $x \in \mathcal{R}$ and $\gamma \in \Gamma$.
\\

The authors would like to thank G\"unter Schlichting and Hartmut F\"uhr for many helpful
discussions.

\section{\large{Earlier Work on Wavelets on Stratified Groups}}

Our results for stratified groups should be contrasted with those of Lemari\'e 
(\cite{Lemarie1}, \cite{Lemarie2}).  He restricted himself to the case where
$\Gamma$ was the set of points all of whose coordinates are integers (to be sure,
this is not always a lattice subgroup).  He constructed an orthonormal basis of
spline wavelets which were $C^N$ (where $N$ is arbitrary, but finite); which had
arbitrarily (but finitely) many derivatives decaying exponentially; and which had
arbitrarily (but finitely) many moments vanishing.  His wavelets were definitely
not smooth; they were built out of splines, that is, functions $\psi$ with $L^M\psi$
a linear combination of Dirac measures for some $M$.  

In this article, we are not seeking orthonormality.   This however enables us
to build in other features which may in certain circumstances be desirable.  
Specifically:
\begin{itemize}
\item
As is well known, the redundancy of a frame is sometimes sought after;
\item
Our continuous wavelets and frames are in the Schwartz space;
\item
In Corollary \ref{C:1} (a), $\phi$ has all moments vanishing and is in the Schwartz space;
\item
In Corollary \ref{C:1} (b), $\phi \in C_c^{\infty}(G)$;
\item
Our prime example, the ``stratified Mexican Hat wavelet'' $Le^{-L/2}\delta$, has the
property that it and all of its derivatives have ``Gaussian'' decay (by the work
of Jersion/ Sanchez-Calle \cite{JerSan} and of Varopoulos \cite{Varapoulos}).
(Here we say a function $F$ on $G$ has ``Gaussian'' decay if for 
some $C, c > 0$,
\[ |F(x)| \leq Ce^{-c|x|^2}. \]
Here $|x|$ is the homogeneous norm of $x$; see section 3 below for homogeneous norms.)
\end{itemize}
 
There are other previous results in wavelet theory on stratified Lie groups,
but -- except in the aforementioned results of Lemari\'e --
high degrees of smoothness and decay, for continuous wavelets
or nearly tight frames, were not previously obtained.  
The existence of admissible functions in $L^2$ was proved by Liu-Peng 
\cite{LiuPeng97} for the Heisenberg group, and by F\"uhr \cite{Fuehr02},
(Corollary 5.28) for general homogeneous groups.  (In contrast to those
works, this article uses no representation theory whatsoever.)
Frames consisting of $L^2$
functions were produced for the Heisenberg group in Maggioni \cite{Maggioni05}.

In the latter article, Maggioni works on a space of homogeneous type which possesses
an involution, and appropriate ``dilations'' and ``translations''; 
examples are the stratified groups considered here (and hypergroups as well).
He assumes that there is an admissible function and creates a wavelet frame
from it.  In the Heisenberg group situation, in order to get an admissible
function, he cites the aforementioned result of Liu-Peng. 
If one instead uses our Theorem \ref{T:2} and Corollary \ref{C:1}, together
with Maggioni's results, one immediately obtains wavelet frames, in the
Schwartz space, on general stratified groups.  One even obtains frames
with the properties stated in our Corollary \ref{C:1} (a) or (b).

In this article, we prefer not to invoke the results of Maggioni, for the
following reason.  Maggioni requires that both the translation parameter 
($b$ in our Theorem  \ref{T:3}) be sufficiently close to $0$ and that the dilation 
parameter $a$ be sufficiently close to $1$.  In Theorem \ref{T:3} (a) we do not need to require
that $a$ be close to $1$; for frames, all that is needed is that Daubechies' criterion 
be satisfied.  This will then enable us to also demonstrate the nearly quadratic convergence
as $a \rightarrow 1$ in Theorem \ref{T:3} (b).

Let us clarify the similarities and differences between our methods and
those of Maggioni, as well as those of earlier authors.
Our method of constructing frames will be through discretizing a continuous problem.
This idea goes back to the beginnings of wavelet theory, for instance, \cite{DGM86} 
and \cite{FJW91}.  In these and other early works, one 
obtained various exact discretizations, where there was no error to be estimated in
replacing an integral by a sum.  More recently, such errors have been estimated, 
specifically in the work of Feichtinger and Gr\"ochenig (\cite{FeichGro1}, \cite{FeichGro2},
\cite{Grochenig91}), 
Gilbert-Han-Hogan-Lakey-Weiland-Weiss (\cite{GHLLWW02}), and
Maggioni \cite{Maggioni05}.  In the latter two references, the error is proved to
have small norm on $L^2$, by use of the $T(1)$ theorem.  
In all of these references, the authors require that both the translation parameter 
($b$ in our Theorem \ref{T:3}) be sufficiently close to $0$ and that the dilation 
parameter $a$ be sufficiently close to $1$.  

We also will use the $T(1)$ theorem.
The reason that we do not have to demand that $a$ be close to $1$, in order
to get a frame, is because we shall discretize, not a continuous
wavelet transform (as in the earlier works just cited),
but rather the operator $R_{\psi}$ which is the operator of convolution with 
$\sum_{j \in \ZZ} \tilde{\psi}_{a^j} * \psi_{a^j}$ (here $\psi = Lf(L)\delta$).
We use the spectral theorem to show that $R_{\psi}$ is bounded below if 
$\xi f(\xi)$ satisfies Daubechies' criterion.   

In section 8 we shall examine wavelet frame expansions in 
other Banach spaces (besides $L^2$).  Again such questions have been discussed in 
the earlier works we have cited (\cite{FeichGro1}, \cite{FeichGro2}, \cite{Grochenig91},
\cite{GHLLWW02}, and 
\cite{Maggioni05}) where again one requires $a$ to be close to $1$.
(In particular, in \cite{Maggioni05}, 
Maggioni addresses such questions on stratified groups.)  Here however we 
shall again require only that the Daubechies criterion be satisfied (so that
$a$ need not be close to $1$).  The novel feature here will be the use of spectral 
multiplier theory (as in \cite{FollandStein82}) to invert $R_{\psi}$ on 
appropriate Banach spaces (such as $L^p$ ($1 < p < \infty$) and the Hardy space
$H^1$).

We also call attention to the important work of Han (\cite{Han00}), on general 
spaces of homogeneous type.  In Theorem 3.35 of 
that article, Han obtains frames by discretizing
a discrete version of the Calder\'on reproducing formula in this general setting.  
He also uses a version of the $T(1)$ theorem to estimate errors.  He
also studies expansions in $L^p$ ($1 < p < \infty$).
However, one cannot expect to obtain nearly tight frames by the
methods in that article.\\

Since we hope this article will be of interest to both the ``wavelet community''
and the ``stratified group community'', we have supplied more details and introductory
material than would be customary had we been writing for only one of these 
communities.\\

In future articles, we will study decay and regularity of dual frames, characterizations
of various Banach spaces through wavelet frame expansion, and analogues of
time-frequency localization for frames.

\section{\large{Notation}}
\label{sec:Notation}

Following \cite{FollandStein82} (which we refer to for further details), we
call a Lie group ${G}$ stratified if it is nilpotent, connected and simply connected, 
and its Lie algebra $\mathfrak{g}$
admits a vector space decomposition $\mathfrak{g}= V_1\oplus \cdots \oplus V_m $ such that
$\left[V_1,V_k\right]=V_{k+1}$ for $1\leq k<m$ and $\left[ V_1,V_m\right]=\{0\}$. \\

If ${G}$ is stratified, its Lie algebra admits a canonical family of dilations, namely
\begin{align}
\delta_r(X_1+X_2+\cdots+ X_m)= rX_1+r^2X_2+\cdots +r^mX_m \quad (X_j\in V_j).\notag
\end{align}

We identify $G$ with $\mathfrak{g}$ through the exponential map.  $G$ is a Lie
group with underlying manifold $\RR^n$, for some $n$.
$G$ inherits dilations from $g$: if $x\in G$ and $r>0$ we write
\begin{equation}\label{eq:dil}
rx=(r^{d_1}x_1,\cdots,r^{d_n}x_n).
\end{equation}
(Here $d_1\leq \cdots\leq d_n$ are those numbers for with
$1\leq k\leq m$ for which $V_k\not=0$). The map $x\rightarrow rx$ is an automorphism of $G$.  \\

The (element of) left (or right) Haar measure on $G$ is simply $dx_1 \ldots dx_n$.  The 
inverse of any $x \in G$ is simply $-x$.  The group law must have the form

\begin{equation} \label{gplaw}
xy = (p_1(x,y),\ldots,p_n(x,y))
\end{equation}
for certain polynomials $p_1,\ldots,p_n$ in $x_1,\ldots,x_n,y_1,\ldots,y_n$.

We let $\mathcal{S}(G)$ denote the space of Schwartz functions on $G$. By definition
$\mathcal{S}(G) = \mathcal{S}(\RR^n)$. \\

The number $Q= \sum_1^m j(dimV_j)$ will be called the \textsl{ homogeneous dimension} of $G$. If
$\phi$ is a function on ${G}$ and $r>0$ , we define $\phi_r$ by
\begin{align}
\phi_r(x)= r^{-Q}\phi(r^{-1}x).
\end{align}

We fix a homogeneous norm function $|\ |$ on $G$ which is smooth away from $0$.  
Thus (\cite{FollandStein82}) $|rx| = r|x|$
for all $x \in G$, $r \geq 0$, $|x^{-1}| = |x|$ for all $x \in G$, and
$|x| > 0$ if $x \neq 0$.  Moreover, for any $a > 0$,
there is a finite $C_a > 0$ such that
$\int_{|x| > R} |x|^{-Q-a} = C_a R^{-a}$ for all $R > 0$.

Let $X_1,\cdots ,X_k$ be a basis for $V_1$ (viewed as left-invariant vector fields on $G$), let
$L= -\sum_1^k X_i^2$
be the sub-Laplacian. This operator (which is hypoelliptic by H\"ormander's theorem
\cite{Hoermander67}) is well known to play on $G$
much the same fundamental role on ${G}$ as (minus) the ordinary Laplacian $\sum_1^N
(\partial_{X_j})^2$ does on $\RR^N$. \\

The operator $L$, restricted to $C_c^{\infty}$, is formally self-adjoint (see Proposition
\ref{zeroint} below).  Its closure has domain 
$\mathcal{D} = \{f \in L^2(G): Lf \in L^2(G)\}$, where
here we take $Lf$ in the sense of distributions.  (This is easily seen through
use of subelliptic estimates.)  From this fact it quickly follows that this closure
is self-adjoint and is in fact the unique self-adjoint extension of $L|_{C_c^{\infty}}$.
We now let $L$ denote this self-adjoint operator.  Suppose that $L$ has spectral resolution

\begin{align}
L= \int_0^\infty \lambda dP_{\lambda}
\end{align}

One then has that $P_{\{0\}} \mathcal{H} = 0$.  To see this, say $f \in L^2(G)$ and
$Lf = 0$; we need to show that $f = 0$.  Since $L$ is the self-adjoint extension
of $L|_{C_c^{\infty}}$, and $Lf = 0$, clearly $Lf = 0$ in the sense of distributions.
But by \cite{GellerLiou}, if $f \in \mathcal{S}'$ and $Lf = 0$, then $f$ is a 
polynomial.  If $f \in L^2(G)$, then surely $f = 0$, as claimed.\\

As usual, if $f$ is a bounded Borel function on $[0,\infty)$, we define the operator $f(L)$ by

\begin{align}
f(L)=\int_0^\infty f(\lambda)dP_{\lambda};
\end{align}

this is well defined and bounded on $L^2({G})$ by the spectral theorem. 
We denote by $f(L)\delta$ the corresponding distribution kernel of the 
bounded operator $f(L)$. Thus

\begin{align}
f(L)\eta= \eta\ast f(L)\delta\quad \forall \;\eta\in \mathcal{S}(G).
\end{align}

\textbf{Notation:} We adopt the $f(L)\delta$ notation, because formally

\begin{align}\label{notation}
f(L)\eta= f(L)\left[\eta\ast\delta\right]= \eta \ast f(L)\delta
\end{align}

since $L$ is left-invariant.\\ \ \\

Let $\RR^+=[0,\infty)$ and set

\begin{align}\notag
\mathcal{S}(\RR^+)=
\{ f\in C^\infty(\RR^+)\;: \forall l, f^{(l)} \;\text{decays rapidly at infinity and}\;\;
lim_{\lambda\rightarrow 0^+} f^{(l)}(\lambda) \;\text{exists}\}.\notag
\end{align}

Then by Borel's theorem on the existence of smooth functions with arbitrary Maclaurin series we
have
$\mathcal{S}(\RR^+)=\mathcal{S}(\RR)|_{\RR^+}$.\\

By \cite{Hulanicki84} (or \cite{Geller80} if G is the Heisenberg group), one has:

\begin{theorem}\label{T:4}
Let $f\in \mathcal{S}(\RR^+)$. Then the distribution kernel of the operator $f(L)=\int_0^\infty
f(\lambda) dP_\lambda$ which we shall denote by $f(L)\delta$ , is a Schwartz function on ${G}$.
\end{theorem}

We have the following elementary lemma on distribution kernels:

\begin{lemma}\label{l:1th} Say $f, g \in \mathcal{S}(\RR^+)$. Then

\begin{enumerate}
\item $\bar{f}(L)\delta= \widetilde{f(L)\delta}$\\
\item $\left[ fg\right](L)\delta= f(L)\delta\ast g(L)\delta$\\
\item For $t>0$ if the function $f^t$ is given by
$f^t(\lambda)= f(t\lambda) \;\;\forall\; \lambda\in [0,\infty)$ ,
then $$\left[ f^t(L)\delta\right]=\left[f(L)\delta\right]_{\sqrt{t}}$$
\end{enumerate}

\end{lemma}

\textbf{Proof:} For $\textit{1}$, using the spectral theorem we have
$\bar{f}(L)= f(L)^\ast$, hence for any
$\phi,\psi\in \mathcal{S}({G})$ we obtain

\begin{align}
<\phi\ast \bar{f}(L)\delta,\psi>&=<\bar{f}(L)\phi, \psi>
=<\phi, f (L)\psi>\\
&=<\phi, \psi\ast f (L)\delta>=
<\phi\ast\widetilde{f(L)\delta},\psi>
\end{align}

which implies the assertion .\\

For $\textit{2}$,
say $\phi \in \mathcal{S}(G)$. By the spectral theorem,

\begin{align}
[(fg)(L)]\phi = g(L)f(L)\phi = [\phi\ast f(L)\delta]\ast g(L)\delta,
\end{align}

yielding $\textit{2}$.\\

For the proof of $\textit{3}$ see Lemma 6.29 of \cite{FollandStein82}.\\

$C$ will always denote a constant, which may change from one occurence to the next.

\section{\large{Proof of Theorem 1.1 and Corollary 1.2}}
\label{sec:ProofOfTheorems11And12}

To prove Theorem ~\ref{T:2}, we need the following lemma:

\begin{lemma}\label{l:2th}

For any $f\in \mathcal{S}(\RR^+)$
with $\int_0^\infty f(s)ds\not=0$
we have
$$K= \int_0^\infty (Lf(L)\delta)_t dt/t=\frac{1}{2}c\delta,$$
where $c= \int_0^\infty f(s)ds$ is a nonzero constant.
\end{lemma}

Note that Theorem ~\ref{T:2} follows immediately from this lemma, since
\begin{align}
(\widetilde{Lf(L)\delta})_t \ast (Lf(L)\delta)_t= \left[\widetilde{Lf(L)\delta} \ast (Lf(L)\delta)
\right]_t,
\end{align}
and by Lemma ~\ref{l:1th} $\widetilde{Lf(L)\delta} \ast (Lf(L)\delta) = Lg(L)\delta $, where
$g(\lambda)= \lambda \mid f(\lambda)\mid^2$. \\

\textbf{Proof 
}
Let $h(\lambda) = \lambda f(\lambda)$.  Write $\psi = h(L) \delta = L f(L)\delta$;
by Lemma \ref{l:1th}, $\psi_t = h^{t^2}(L) \delta$ for any $t > 0$.
Define
$K_{\epsilon , A}= \int_\epsilon^A \psi_t dt/t$. Since $\int_G \psi =
\int_{{G}} Lf(L)\delta =0$,
by Theorem 1.65 \cite{FollandStein82}, $ \int_\epsilon^A \psi_t dt/t$ converges in
$\mathcal{S}'$ as $\epsilon \rightarrow 0$ and $A\rightarrow \infty$ to the tempered distribution
$K= \int_0^\infty \psi_t dt/t$, which is $C^\infty$ away from $0$.
Suppose $\phi_1 \in \mathcal{S}(G)$. Then $\phi_1\ast K_{\epsilon , A}\in \mathcal{S}$ and for any
$\phi_2\in \mathcal{S}(G)$ we have

\begin{align} \notag
<\phi_1\ast K_{\epsilon , A},\phi_2>= <K_{\epsilon , A}, \widetilde{\phi}_1\ast\phi_2>
&= \int_\epsilon^A <\psi_t ,\widetilde{\phi}_1\ast\phi_2 >dt/t\\ \notag
&=\int_\epsilon^A <\phi_1\ast \psi_t ,\phi_2>dt/t\\ \notag
&=\int_\epsilon^A <[h^{t^2}(L)]\phi_1 ,\phi_2>dt/t\\ \notag
&=\int_\epsilon^A \int_0^{\infty} t^2 \lambda f(t^2 \lambda) d<P_{\lambda}\phi_1,\phi_2>dt/t\\ \notag
&=\int_0^{\infty} \int_\epsilon^A t^2 \lambda f(t^2 \lambda) dt/t\: d<P_{\lambda}\phi_1,\phi_2>\\ \notag
&=\frac{1}{2}\int_0^{\infty}\int_{\lambda\epsilon^2}^{\lambda A^2} f(t)dt\:
d<P_{\lambda}\phi_1, \phi_2>.\notag
\end{align}

Letting $F(x)= -\int_x ^\infty f(s)ds$ (so that $F'=f$) we see that this double integral equals

\begin{align}
\int_0^{\infty} \int_{\lambda\epsilon^2}^{\lambda A^2} f(t)dt\: d<P_{\lambda}\phi_1,\phi_2>
= \int_0^{\infty} \Big( F(\lambda A^2)-F(\lambda\epsilon^2) \Big)d<P_{\lambda}\phi_1,\phi_2>.
\end{align}

Since the function $F$ is bounded, and the measure $<P_{\lambda}\phi_1,\phi_2>$
is supported on $(0,\infty)$ (in that $P_{\{0\}}= 0$), we see that

\begin{align}
lim_{\epsilon\rightarrow 0\;A\rightarrow \infty} \int_0^{\infty} \Big(F(\lambda
A^2)-F(\lambda\epsilon^2)\Big)d<P_{\lambda}\phi_1,\phi_2> &=
\int_0^{\infty} \int_0^\infty f(s)ds\: d<P_{\lambda}\phi_1,\phi_2>\\
&= \int_0^\infty f(s)ds <\phi_1,\phi_2>.
\end{align}

This proves the Lemma.  Thus Theorem ~\ref{T:2} is established as well.\\
\ \\
To begin the proof of Corollary \ref{C:1}, if $\alpha = (\alpha_1,\ldots,\alpha_n)$ is a
multi-index, we
let $|\alpha| = \sum_k d_k \alpha_k$. Note that $|\alpha|$ is the homogeneous degree of the
monomial
$x^{\alpha}$, since $(rx)^{\alpha} = r^{|\alpha|}x^{\alpha}$ for $r > 0$.
 For any positive integer $k$, $L^k x^{\alpha}$ is a polynomial which is homogeneous of
degree $|\alpha| - 2k$; it must therefore be identically zero if $|\alpha| - 2k < 0$.
Integration by parts now at once shows the following proposition:

\begin{proposition}
\label{P:0}
If $F \in \mathcal{S}(G)$, and if $|\alpha| < 2k$, then $\int_G x^{\alpha} L^k F = 0$.
\end{proposition}

{\bf Proof of Corollary \ref{C:1}} For (a), select any nonzero $g \in \mathcal{S}(\RR^+)$ which
vanishes
identically in a neighborhood of $0$. For any positive integer $k$, define $g_k(x) = g(x)/x^k$;
then $g_k \in \mathcal{S}(\RR^+)$, and $g(L)\delta = L^k g_k(L) \delta$. By Theorem \ref{T:2} and
Proposition
\ref{P:0}, $g(L) \delta$ is admissible and has all moments vanishing. \\

For (b), we note that if $g \in C_c^{\infty}(\RR)$ is real-valued and even, and if
$m(\lambda) = \hat{g}(\sqrt{\lambda})$, then $m(L)\delta \in C_c^{\infty}(G)$. (This is proved in
the appendix to \cite{GrafakosLi00}; the argument is there attributed to J. Dziubanski, but he says
the
result was well-known; it appears to be based on ideas of Michael Taylor.)
Thus, if $g \neq 0$, then for any positive integer $k$, $ \phi_k = L^k
m(L)\delta
= L(L^{k-1} m(L) \delta)$ is admissible and in
$C_c^{\infty}(G)$, and $\int x^{\alpha} \phi_k = 0$ whenever $|\alpha| < 2k$. (Note that $\phi_k$
cannot be
identically zero, for then $\lambda^k m(\lambda)$ would be identically zero, so $g$ would be
zero.) This
completes the proof.\\
\ \\
{\bf Remark} Corollary \ref{C:1} (b) improves on Theorems 1.61 and 1.62 of Folland-Stein
\cite{FollandStein82},
at least for stratified $G$. There it was shown that there exist
$\phi^1,\ldots,\phi^M,\psi^1,\ldots,\psi^M
\in \mathcal{S}(G)$ with arbitrarily many moments vanishing, with the $\psi^j$ having compact
support, and
with $\sum_1^M \int_0^{\infty}\phi_t^j * \psi_t^j dt/t = \delta$; here $M$ depended on the number
of moments
one wanted to vanish. Now we see that we can always take $M=1$ and $\psi^1 = \tilde{\phi^1}$, so
that
both have compact support.

\section{\large{Continuous Wavelet Transform}}
\label{sec:ContinuousWaveletTransform}

In this section we study the continuous wavelet transform with respect to the quasiregular
representation
of the group $M:={G} \ltimes (0,\infty )$ , where ${G}$ is a stratified group with homogeneous
degree $Q$ and with Haar measure $db$. \\
$M$ is a locally compact group with left Haar measure $d\mu(M)=a^{-(Q+1)}dadb $. \\

The positive number $a$ defines an automorphism of the group $G$, which acts by dilation. The
quasi-regular representation $\pi $ of $M$ acts on
$L^2({G})$ as follows:\\

Let $\phi\in L^2(G)$, then

\begin{align}\label{q:44}
(\pi (x,a) \phi )(y)= (T_xD_a\phi )(y)=a^{-Q/2}\phi(a^{-1}(x^{-1}y))\;\;\forall x,y\in {G}
\;,\; \forall a>0
\end{align}

Thus $T_x$ acts by left translation by $x^{-1}$ , while $D_a$ denotes a unitary dilation operator
with respect to $a$.\\

The following definition and more details can be found for example in \cite{Fuehr02}.
\begin{definition}\label{de:fuehr}
Let $\phi$ and $\psi$ be any fixed functions in $L^2(G)$. 
Define the coefficient function
$V_{\phi,\psi}$ on $G$ by
\begin{align}
V_{\phi,\psi}: (x,a)\mapsto <\psi,T_xD_a\phi >
\end{align}

The coefficient function $V_{\phi,\psi}$ is not necessarily square integrable on $M$.
The function $\phi$ is called admissible when for any $\psi$ the associated
coefficient function $V_{\phi,\psi}$ is square integrable, and the operator
\begin{align}V_\phi: L^2 (G) \longrightarrow L^2(M),
\end{align}

given by $[V_{\phi}(\psi)](x,a) = V_{\phi,\psi}(x,a)$,
is an isometry.
Then, for the admissible vector $\phi$, the bounded operator $V_{\phi}$ is 
called a continuous wavelet transform of $L^2(G)$.
\end{definition}

We shall soon show (in Proposition ~\ref{pro:1th} below) 
that this (accepted) definition of
admissible is consistent with our usage of the word \textit{admissible} in Theorem ~\ref{T:2}.

The existence of admissible vectors in $L^2( G)$ for   $\pi$  was proved by F\"uhr \cite{Fuehr02},(
Corollary 5.28) for homogeneous groups. We recall this in the next Theorem:

\begin{theorem}\label{Fr}
Let $M=G\ltimes H$, where $G$ is a homogeneous Lie group and $H$ is a one-parameter group of
dilations. Then the quasi-regular representation $\pi$ is contained in the left regular
representation $\lambda_M$. Hence there exists a continuous wavelet transform on $G$ arising from
the action of $G$ by left translations and the action of the dilations.
\end{theorem}

We now show 
(without use of Theorem \ref{Fr})
that there exist admissible $\phi \in
\mathcal{S}(G)$. We claim:

\begin{proposition}\label{pro:1th}
Say $\phi\in \mathcal{S}(G)$ and $\int \phi=0$, so that by Theorem 1.65 of \cite{FollandStein82},
if
\begin{align}
K_{\epsilon,A}= \int_\epsilon^A\frac{\tilde{\phi_t} \ast \phi_t}{t}dt 
\end{align}

then $K= lim_{\epsilon\rightarrow 0, A\rightarrow \infty}K_{\epsilon,A}$ exists in
$\mathcal{S}'(G)$ , $C^\infty$ away from $0$ and is homogeneous of degree $-Q$. Then $\phi$ is
admissible (in the sense of Definition ~\ref{de:fuehr} ) if and only if $K=\delta$ up to a
constant multiple. In particular if $0\not= f\in \mathcal{S}(\RR^+)$, then $\phi=Lf(L)\delta$ is
admissible.
\end{proposition}

\textbf{Proof:} For $\psi\in L^2({G})$ we have:

\begin{align}\label{q:46}
\int_M \mid V_{\phi} \psi\mid^2 &
=\int_{{G}}\int_0^\infty \mid <\psi,T_bD_a\tilde{\phi}>\mid ^2 d\mu (M)\\ \notag
&=\int_{{G}}\int_0^\infty \mid \psi\ast (D_a\tilde{\phi})(b) \mid ^2 a^{-(Q+1)}dadb \quad {\rm
so},\\
\label{eq:*}
\int_M \mid V_{\phi} \psi\mid^2&=
\int_0^\infty
\parallel \psi\ast (D_a\tilde{\phi} )\| ^2a^{-(Q+1)}da.
\end{align}

But for any $a>0$,

\begin{align}\label{eq:**}
\parallel \psi\ast (D_a\tilde{\phi})\| ^2 a^{-Q}
= <\psi,\psi\ast \widetilde {(D_a\phi)} \ast (D_a\phi) >
a^{-Q}= <\psi,\psi\ast (\widetilde{\phi_a} \ast \phi_a)>.
\end{align}

Since $K_{\epsilon,A} \rightarrow K$ in $\mathcal{S}'$, if $g\in \mathcal{S}$, then $g\ast
K_{\epsilon,A} \rightarrow g\ast K$ pointwise and for some $N,C$

\begin{align}\label{eq:***}
\mid (g \ast K_{\epsilon,A})(x)\mid \leq C(1+\mid x\mid)^N \quad \quad \text{for all } x,\epsilon,
A.
\end{align}

Using the dominated convergence theorem in (~\ref{eq:*}) and (~\ref{eq:**}) , if $\psi \in
\mathcal{S}(G)$,then

\begin{align}\label{eq:psiK}
\parallel V_\phi \psi\parallel_{L^2}^2 = <\psi, \psi\ast K>\leq C\parallel \psi\parallel_{L^2}^2
\end{align}

since the map $\psi\rightarrow \psi\ast K$ is bounded on $L^2(G)$. $V_\phi$ thus maps
$\mathcal{S}(G)$
to $L^2(M)$ and has a unique bounded extension to a map from $L^2(G)$ to $L^2(M)$. But if
$\psi_k\rightarrow \psi$ in
$L^2(G)$, surely $V_\phi \psi_k\rightarrow V_\phi \psi$ pointwise , so this extension can be none
other than $V_\phi$. Accordingly (~\ref{eq:psiK}) holds for all $\psi\in L^2(G)$.  
We thus have

\begin{align}
\parallel V_\phi \psi\parallel_{L^2(M)}= \parallel \psi\parallel_{L^2} \forall \psi\in L^2
&\Longleftrightarrow <\psi, \psi\ast K> = <\psi,\psi> \quad \forall \psi\in L^2 \\
&\Longleftrightarrow \psi\ast K=\psi \quad \forall \psi\in L^2 \\
&\Longleftrightarrow K=\delta \quad \text{up to a constant.}\quad
\end{align}

as desired.  (In the second implication, we have used polarization.)  
This completes the proof. \\

\section{Lemmas on Vector Fields}

In this section we gather a number of facts which will be needed in our discussion
of frames.  These facts are analogues for $G$ of very standard facts on $\RR^n$ 
(such as the fundamental theorem of calculus -- see Lemma \ref{ftc} below).

The right-invariant vector fields $Y_l$ ($1 \leq l \leq n$) may be defined by
\[ Y_l g = -\widetilde{X_l\tilde{g}}  \]
for $g \in C^1(G)$.  

We note:

\begin{proposition}\label{zeroint}
Suppose $\phi \in \mathcal{S}(G)$.  Then, for all $l$, $\int_G X_l \phi = 0$
and $\int_G Y_l \phi = 0$.
\end{proposition}
{\bf Proof} Note that each $X_l$ is homogeneous of degree $a_l$.  
This forces $X_l$ to have the form
\[ X_l = \frac{\partial}{\partial x_l} + 
\sum_{k > l}p_k(x) \frac{\partial}{\partial x_k},  \]
where $p_k$ is a homogeneous polynomial of degree $a_k - a_l < a_k$.  
(See \cite{FollandStein82} for a detailed proof of this.)  Accordingly
$p_k(x)$ must actually be a polynomial in $x_1,\ldots,x_l$, so multiplication by
it commutes with $\partial/\partial x_k$ for $k > l$.

Accordingly $\int_G X_l \phi = 0$. By using $\tilde{\ }$ we see that
$\int_G Y_l \phi = 0$ as well.\\

If $x = (x_1,\ldots,x_n) \in G$, and
$t > 0$, for want of a better notation, let us define
\[  [t]x = (tx_1,\ldots,tx_n)  \]
(recall that $tx$ means something else, see (\ref{eq:dil})).  

Recall that we are identifying $G$ with 
$\mathfrak{g}$ through the exponential map.  Then, if $x
\in G$, we say that the point $\exp(x \cdot X)(0)$ has coordinates $x$.

\begin{lemma} \label{ftc}
(a) Suppose that $x \in G$ and that $U$ is an open neigborhood
of the line segment \newline
$\{[t]x: 0 \leq t \leq 1\}$.  
If $g \in C^1(U)$, then 
\begin{equation}\label{eq:52L}
g(x)-g(0) = \int_0^1 [(x \cdot X) g]([t]x) dt
\end{equation}
and
\begin{equation}\label{eq:52}
g(x)-g(0) = \int_0^1 [(x \cdot Y) g]([t]x) dt.
\end{equation}

(b) Suppose that $x, u \in G$ and that $U$ is an open neighborhood
of the set $\{u([t]x): 0 \leq t \leq 1\}$.
If $h \in C^1(U)$, then
\begin{equation}\label{eq:52LL}
h(ux)-h(u) = \int_0^1 [(x \cdot X) h](u([t]x)) dt.
\end{equation}
(c) Suppose that $x, u \in G$ and that $U$ is an open neighborhood
of the set $\{([t]x)u: 0 \leq t \leq 1\}$.
If $h \in C^1(U)$, then
\begin{equation} 
\label{eq:52RR}
h(xu)-h(u) = \int_0^1 [(x \cdot Y) h](([t]x)u) dt.
\end{equation}
\end{lemma}
{\bf Proof} For (a), we note that
\begin{eqnarray*}
g(x) - g(0) &=& g(\exp(x \cdot X)(0)) - g(0)\\
&=& \int_0^1 \frac{d}{dt}g(\exp(t[x \cdot X])(0)) dt\\
&=& \int_0^1 [(x \cdot X) g]([t]x) dt.
\end{eqnarray*}
proving (\ref{eq:52L}).  Applying $\tilde{\ }$ to (\ref{eq:52L}), we find
(\ref{eq:52}) as well.  For (b), we apply (\ref{eq:52L}) to the function $g = h_u$
where $h_u(x) = h(ux)$.  For (c) we apply (\ref{eq:52}) to the function $g = h_u$
where $h_u(x) = h(xu)$.  This completes the proof.

We will be needing two applications of Lemma \ref{ftc}, Propositions \ref{nbfmvt}
and \ref{czmvt} below.  First, however, some remarks on homogeneous norms.

A homogeneous norm function satisfies a type of triangle
inequality (\cite{FollandStein82}, equation (1.8)):
for some $C > 0$, $|xy| \leq C(|x|+|y|)$ for all $x,y \in G$.   
We shall need three other facts about homogeneous norms:

\begin{proposition} \label{triget}
There exists $c > 0$  such that for all $R > 0$, if
$|u^{-1}x| \geq 2R$, then 
\newline
$\min_{|u^{-1}y| \leq R} |x^{-1}y|
\geq c|u^{-1}x|$.    
\end{proposition}
{\bf Proof}  Since $x^{-1}y = (x^{-1}u)(u^{-1}y) 
= (u^{-1}x)^{-1}(u^{-1}y)$,
we may after a translation assume $u = 0$.  
It is enough to show that, for some $c > 0$, if $|y| \leq |x|/2$, then
$|x^{-1}y| > c|x|$.  After a dilation we may assume
$|x|=2$ and $|y| \leq 1$.  By the triangle inequality, for some $C>0$,
$|x^{-1}y| \geq |x|/C -|y|$, so we may assume also that $|x| \leq 2C$.
But $|x^{-1}y|$ does not vanish for $(x,y)$ in the compact set
$\{x: 2 \leq |x| \leq 2C\} \times \{y: |y| \leq 1\}$, so it has
a positive minimum there, as desired.\\

Sometimes we use the ``standard homogeneous norm function'' on $G$, defined by
\[ |x| = (\sum_{k=1}^n |x_k|^{2b_k})^{1/2A},  \]
where $A = a_1 \ldots a_n$, and each $b_k = A/a_k$.  We shall clearly 
indicate when we do this.

\begin{proposition} \label{hmnm1}
There is a constant $C > 0$ such that for all $x = (x_1,\ldots,x_n) \in G$,
$|x_m| \leq C|x|^{a_m}$ for $1 \leq m \leq n$.
\end{proposition}

\begin{proposition} \label{hmnm2}
There is a constant $C > 0$ such that for all $x \in G$ and all
$t$ with $0 \leq t \leq 1$, we have $|[t]x| \leq C|x|$.  If $|\ |$ is
the standard homogeneous norm function, we can take $C = 1$.
\end{proposition}
{\bf Proof of Propositions \ref{hmnm1} and \ref{hmnm2}}  Since any two 
homogeneous norms are equivalent, we may assume that $|\ |$ is the
standard homogeneous norm function.  But in that case the propositions are evident 
(and we can take $C=1$ in both).\\

We now turn to the applications of Lemma \ref{ftc}.  We define
a {\em normalized bump function} to be a $C^1$ function with support in the
unit ball  $B(0,1) =
\{x: |x| < 1\}$ with $C^1$ norm less than or equal to $1$. For  any
function $f: G \rightarrow \CC$,
if $R > 0$ and $u \in G$, we let  $f^{R,u}(x) = f(R^{-1}(u^{-1}x))$.
We claim:

\begin{lemma}  \label{nbfmvt}
There exists a constant $C > 0$ such that for all normalized bump functions
$f$, all $R > 0$, and all $u, x, y \in G$ we have
\[ |f^{R,u}(xy) - f^{R,u}(x)|  \leq C\sum_{k=1}^n \frac{|y_k|}{R^{a_k}}.\]
\end{lemma}
{\bf Proof} We have
\[ f^{R,u}(xy) - f^{R,u}(x) = f(R^{-1}(u^{-1}xy)) - f(R^{-1}(u^{-1}x)) 
= f(x^{\prime} y^{\prime}) - f(x^{\prime}),  \]
where $x^{\prime} = R^{-1}(u^{-1}x)$ and 
\[ y^{\prime} = R^{-1}y = (\frac{y_1}{R^{a_1}},\ldots, \frac{y_n}{R^{a_n}}). \]
In proving the lemma we may therefore assume that $R=1$ and $u = 0$, so 
that $f^{R,u} = f$.  In that case we use Lemma \ref{ftc} (b) to find that
\[ |f(xy)-f(x)| = \int_0^1 \left[(y \cdot X) f\right](x([t]y)) dt \leq C\sum_{k=1}^n |y_k|  \]
as claimed, since the functions $X_k f$ are bounded (uniformly for all
normalized bump functions $f$).\\

We now turn to our second application of Lemma \ref{ftc}.  First we
define a {\em Calderon-Zygmund kernel} to be a complex-valued function
$K(x,y)$, defined for all $x, y \in G$ with $x \neq y$, which is continuous
(off the diagonal), and which, for some $C,c > 0$, satisfies the following 
three estimates (for all $x,y \in G$ with $x \neq y$):

\begin{equation} \label{cz1}
|K(x,y)| \leq \frac{C}{|y^{-1}x|^{Q}};
\end{equation}

\begin{equation} \label{cz2}
\mbox{If } |x^{-1}x'| \leq c|y^{-1}x|, \mbox{ then } |K(x',y) - K(x,y)|
\leq C\frac{|x^{-1}x'|}{|y^{-1}x|^{Q+1}};
\end{equation}

\begin{equation} \label{cz3}
\mbox{If } |y^{-1}y'| \leq c|y^{-1}x|, \mbox{ then } |K(x,y') - K(x,y)| 
\leq C\frac{|y^{-1}y'|}{|y^{-1}x|^{Q+1}}.
\end{equation}

We then claim:

\begin{proposition} \label{czmvt}
Suppose $K(x,y)$ is defined and $C^1$ away from the diagonal in $G \times G$,
and that for some $A > 0$, 
\begin{equation} \label{czc1}
|X_x^{\alpha} X_y^{\beta}K(x,y)| \leq
A|y^{-1}x|^{-Q-|\alpha|-|\beta|}.
\end{equation}
whenever $0 \leq \alpha_1 + \ldots + \alpha_n + \beta_1 + \ldots + \beta_n \leq 1$,
and whenever $x, y \in G$ with $x \neq y$.
Then $K$ is a Calderon-Zygmund kernel.
(Here
$X_x^{\alpha} = X_1^{\alpha_1}\ldots X_n^{\alpha_n}$, where the $X_k$ are taken
in the $x$ variable.)\\
\end{proposition}
{\bf Proof} By taking $\alpha = \beta = 0$ in (\ref{czc1}), we have (\ref{cz1}).
To prove (\ref{cz2}), we may assume we are using the standard homogeneous norm
function; we will then show that (\ref{cz2}) holds with $c = 1/2$.  

In this proof, it will be convenient let $X_{k,1} K(x,y)$ denote the result of 
applying $X_k$ to $K$ in the $x$ variables.

Suppose $x \neq y$ and $|x^{-1}x'| \leq |y^{-1}x|/2 = |x^{-1}y|/2$.
If $0 \leq t \leq 1$, then
by Proposition \ref{hmnm2}, $|[t](x^{-1}x')| \leq |x^{-1}x'| \leq |x^{-1}y|/2$ 
as well.  In particular, $[t](x^{-1}x') \neq x^{-1}y$, so $x([t](x^{-1}x')) \neq y$.
Moreover, by Proposition \ref{triget}, there exists a $c_1 > 0$ (independent of
the specific values of $x,y,x',t$) such that
\[  |y^{-1}x([t](x^{-1}x'))| \geq c_1|y^{-1}x|.  \]

We write $x' = x(x^{-1}x')$.  
Using Lemma \ref{ftc} and Proposition \ref{hmnm1}, we find that for some $C_1, C_2, C_3 > 0$,

\begin{eqnarray*}
|K(x',y) - K(x,y)| &=& |\int_0^1 \sum_{k=1}^n (x^{-1}x')_k (X_{k,1} K(x([t](x^{-1}x')),y) dt|\\
& \leq & C_1 A\sum_{k=1}^n \frac{|(x^{-1}x')_k|}{|y^{-1}x|^{Q+a_k}}\\
& \leq & C_2 A\sum_{k=1}^n \frac{|x^{-1}x'|^{a_k}}{|y^{-1}x|^{Q+a_k}}\\
& \leq & C_3 A\frac{|x^{-1}x'|}{|y^{-1}x|^{Q+1}}
\end{eqnarray*}
so that (\ref{cz2}) holds.  (Note for later purposes that $C_1, C_2, C_3$ depend only 
on the group $G$ and not in any way on $K$.)  The proof of (\ref{cz3}) is exactly
analogous.  This proves the proposition.\\

We will be using Lemma \ref{nbfmvt} and Proposition \ref{czmvt} in conjunction with 
the $T(1)$ theorem for stratified groups.  We review this theorem in a moment.

First, however, some definitions.  Suppose that a linear operator $T: C_c^{1}(G) 
\rightarrow L^2(G)$.  One says that $T$ is {\em restrictedly
bounded} if there is a $C > 0$ such that $\|T(f^{R,u})\|_2 \leq CR^{Q/2}$ for all
normalized bump functions $f$, all $R$ and all $u$.  

If $T: C_c^{1} \rightarrow L^2(G)$ is linear,
we say that a linear operator $T^* : C_c^{1} \rightarrow L^2(G)$
is its formal adjoint if for all $f, g \in C_c^{1}$ we have
\begin{equation}  \label{formadj}
<Tf,g> = <f,T^*g>.
\end{equation}
$T^*$ is evidently unique if it exists. 

We will be using the ``easier case'' of the David-Journ\'e $T(1)$ theorem 
\cite{DavidJourne84} for stratified groups (\cite{Lemarie84}
or \cite{Stein93}, pages 293-300).  (The latter reference is only for $G=\RR^n$, but the 
proof for general $G$ requires only minor changes -- see the appendix to this paper.) 
We may formulate this theorem as follows:

\begin{theorem} \label{tof1}
Suppose that $T: C_c^{1}(G) \rightarrow L^2(G)$ has a formal adjoint
$T^*: C_c^{1}(G) \rightarrow L^2(G)$.  Suppose further:\\
(i) $T$ and $T^*$ are restrictedly bounded; \\
(ii) There is a Calderon-Zygmund kernel $K$ such that if $f \in C_c^{1}$, 
then for $x$ outside the support of $f$, $(Tf)(x) = \int K(x,y)f(y)dy$; and \\
(iii) $T(1) = T^*(1) = 0$.\\
Then $T$ extends to a bounded operator on $L^2$.
\end{theorem}

Condition (iii) means precisely
(\cite{Stein93}, pages 300-301) that whenever $f \in C_c^{\infty}(G)$
and $\int_G f = 0$, we have that $\int_G Tf = \int_G T^*f = 0$.  In fact
we shall show that this is true for all $f \in C_c^1(G)$ with $\int_G f = 0$.
(Even without condition (iii), conditions (i) and (ii) imply that for all 
such $f$, $Tf$ and $T^*f$ are in $L^1(G)$.)

In fact we shall need a quantitative version of Theorem \ref{tof1}:

\begin{theorem} \label{tof1q}
There exist $C_0, N > 0$, such that for any $A > 0$, we have the following.
Whenever $T: C_c^{1}(G) \rightarrow L^2(G)$ has a formal adjoint
$T^*: C_c^{1}(G) \rightarrow L^2(G)$, and whenever $T, T^*$ satisfy:\\
\ \\
(i) $\|T f^{R,u}\|_2 \leq A R^{Q/2}$ and $\|T^* f^{R,u}\|_2 \leq A R^{Q/2}$ 
for all normalized bump functions $f$;\\
(ii) There is a kernel $K(x,y)$, $C^1$ off the diagonal, 
such that if $f \in C_c^{1}$, 
then for $x$ outside the support of $f$, $(Tf)(x) = \int K(x,y)f(y)dy$; and 
whenever at most one of $\alpha_1, \ldots, \alpha_n, \beta_1, \ldots, \beta_n$
is not zero, and whenever $x, y \in G$ with $x \neq y$, we have
\[ |X_x^{\alpha} X_y^{\beta}K(x,y)| \leq A|y^{-1}x|^{-Q-|\alpha|-|\beta|};\ \ \ \  {\rm and} \]
(iii) $T(1) = T^*(1) = 0$,\\
\ \\
then $T$ extends to a bounded operator on $L^2$, and $\|T\| \leq C_0 A$.

\end{theorem}
{\bf Proof}   This follows at once from an examination of the proofs of Theorem
\ref{tof1} (in \cite{Lemarie84} or \cite{Stein93}) and of Proposition 
\ref{czmvt} above.

\section{Frames}

Suppose now that one has a discrete subset $\Gamma$ of $G$, and a bounded
measurable set $\mathcal{R} \subseteq G$ of positive measure, such that every
 $g \in G$ may be written uniquely
in the form $g = x\gamma$ with $x \in \mathcal{R}$ and $\gamma \in \Gamma$.\\

For example, one could choose $\Gamma$ to be any lattice
subgroup of $G$, if one is available. Thus $\Gamma$ is a discrete subgroup of $G$, such that
$G/\Gamma$ is compact. (Note: by \cite{CorwinGreenleaf89}, page 197, equation (2),
it is equivalent to assume that $\Gamma$ is a discrete
subgroup of $G$, such that $G/\Gamma$ has
finite volume with respect to the induced invariant measure. 
If the coefficients of all the polynomials appearing in
(\ref{gplaw}) are integers, as is the case for the Heisenberg group,
one could take $\Gamma$ to be the integer lattice, namely
the set of points all of whose coordinates are integers.)
We then let $\mathcal{R}$ be a fundamental region
for $G/\Gamma$. (By this we
mean a bounded measurable subset of $G$,  
of positive measure, consisting of precisely one representative of
each right coset of $\Gamma$.) Thus every $g \in G$ may be written uniquely
in the form  $g=\gamma x$ with
$x \in \mathcal{\mathcal{R}}$, $\gamma \in \Gamma$.)

\begin{definition}\label{definition-of-frame}
A countable subset $\{e_n\}_{n\in I}$ of a Hilbert space $\mathcal{H}$ is said to be a
\textbf{frame} if there exist two positive numbers $A\leq B$ such that, for any
  $f\in\mathcal{H}$,
  \begin{align}
 A\parallel f\parallel^2 \leq \sum_{n\in I} \mid <f,e_n>\mid^2 \leq B\parallel f\parallel^2.
  \end{align}
the positive numbers $A$ and $B$ are called  \index{frame! bounds} {\em frame
bounds}.
\end{definition}
Note that the frame bounds are not unique. The \textit{a lower frame bound} is the supremum
over all lower  frame bounds, and the  \textit{optimal  upper  frame  bound} is the infimum over all
upper frame bounds. The optimal frame bounds are of course frame bounds.
 The frame is called a {\em tight frame}  when we can take
$A=B$.  (Informally, we also say the frame is {\em ``nearly tight''} if
$B/A$ is close to $1$.)  Frames were introduced in \cite{DuffinSchaeffer52}.

We consider $\phi \in \mathcal{S}(G)$ with $\int \phi = 0$.
For $a,b > 0$, we define
\[ \phi_{j,b\gamma}(x) = [D_{a^j}T_{b\gamma} \phi](x) = 
a^{-jQ/2} \phi([b\gamma]^{-1}[a^{-j}x]). \]
($a$ will usually be fixed.) The set $\{\phi_{j,b\gamma}\}$
is called the wavelet system generated by $\phi$.  
We seek conditions on $\phi$
and the numbers $a,b > 0$ which guarantee that this wavelet
system is a frame (in which case it is called a wavelet frame).

In order to do this we study the operator

\begin{align}\notag
S_{\phi,b}:\;  f\rightarrow \sum_{\gamma\in\Gamma, j\in\ZZ}
<f, \phi_{j,b\gamma}> \phi_{j,b\gamma}
\end{align}

It is not hard to see \cite{Daubechies92} that
$\{\phi_{j,b\gamma}\}$ is a frame if and only if:
for any $f \in L^2(G)$, 
the series defining $S_{\phi,b}f$ converges unconditionally 
to $f$ in $L^2(G)$; and $S_{\phi,b}$ 
is bounded on $L^2(G)$; and $S_{\phi,b} \geq AI$  for some 
{\em strictly positive} number $A$. (If the frame is ``nearly tight'', that is
if, for certain frame bounds $A,B$ one knows that $B/A - 1 = \epsilon$ is small, then 
(1/A)$S_{\phi,b} f$ is a good approximation to $f$.  Indeed, for any $f \in L^2$,
$A\|f\|^2 \leq <S_{\phi,b}f,f> \leq B\|f\|^2$ implies that, as operators,
$0 \leq (1/A)S_{\phi,b} - I \leq \epsilon I$, whence $\|(1/A)S_{\phi,b} - I\|
\leq \epsilon$.  For this
reason, one generally prefers ``nearly tight'' frames.)

More generally we shall need to consider $\phi, \psi \in \mathcal{S}(G)$
with
$\int \phi = \int \psi = 0$ and look at operators of the form

\begin{align}\notag
S_{\phi,\psi,b}:\; f\rightarrow  \sum_{\gamma\in\Gamma, j\in\ZZ}
<f,\phi_{j,b\gamma}> \psi_{j,b\gamma}.
\end{align}

\begin{theorem}\label{framebdd}
Fix $a > 0$. In parts (a), (b), (c) and (d) we also fix $\phi, \psi \in
\mathcal{S}(G)$ with
$\int \phi = \int \psi =  0$.\\
(a) For any $0 < b < 1$ and $f \in C_c^{1}(G)$, the series  defining
$S_{\phi,\psi,b}f$
converges absolutely, uniformly on $G$.\\

(b) For any $0 < b < 1$ and $f \in C_c^{1}(G)$,  $S_{\phi,\psi,b}f \in
L^2(G)$.\\

(c) For some $C > 0$,  $\|S_{\phi,\psi,b}f\|_2 \leq (C/b^Q)\|f\|_2$ for all
$0 < b < 1$ and $f  \in C_c^{1}(G)$.\\

Consequently $S_{\phi,\psi,b}$ extends to be a bounded  operator on
$L^2(G)$.  (In fact, if we put $T= S_{\phi,\psi,b}$, then
$T$ satisfies the hypotheses of Theorem \ref{tof1}.)\\

(d) If $f,g \in L^2(G)$, then 
\begin{equation}\label{eq:wk}
<S_{\phi,\psi,b}f,g>= 
\sum_{\gamma\in\Gamma, j\in\ZZ}<f, \phi_{j,b\gamma}>
<\psi_{j,b\gamma},g>;
\end{equation}
here the series converges absolutely.\\

(e) Say $\mathcal{B}_0$ is a bounded subset of $\mathcal{S}(G)$, $f \in C_c^{1}(G)$ and $b > 0$.
Then the series defining $[S_{\phi,\psi,b}f](x)$ converges absolutely, uniformly for
$x \in G$ and $\phi, \psi \in \BB_0$ with $\int \psi = \int  \phi = 0$.  \\

(f) If $\mathcal{B}_0$ is a bounded subset of $\mathcal{S}(G)$,  then there exists a
constant $C$
such that $\|S_{\phi,\psi,b}\| \leq C/b^Q$  for all $0 < b < 1$ and all
$\psi, \phi \in \BB_0$ with
$\int \psi = \int  \phi = 0$.
\end{theorem}
{\bf Remark} For the boundedness of $S_{\phi,\psi,b}$ on $L^2$, one may also
consult section 6 of Maggioni \cite{Maggioni05}.   If $G= \RR^n$,
the fact that $T= S_{\phi,\psi,b}$ satisfies the conclusions
of Theorem \ref{tof1} has a long history.  For instance, in Lemma 9.1.5 of 
\cite{Daubechies92}, condition {\em (ii)} of Theorem \ref{tof1} (b) is verified 
for this $T$, if $G = \RR$.  If $G = \RR^n$, Theorem \ref{tof1} is verified for
this $T$, and more general operators $T$, in \cite{GHLLWW02}, sections 2.1-2.3.
\ \\

To prove the theorem, we shall need the following technical lemma.
\begin{lemma}\label{gmlem}
For $N > 0$ define the function $g_N$  on $G$ by
\[ g_N(x) = (1 +|x|)^{-N}. \]
Then:\\
(a) Let $B_0$ be a  bounded subset of $G$.  Then for some $C >0$,
\[ g_N(x) \leq C  g_N(y^{-1}x) \]
for all $x \in G$ and $y \in B_0$.\\
(b) Say $M, N >  Q/2$, and suppose $0 < L < \min(M-Q/2, N-Q/2)$.  Then for
some $C  > 0$,
\[ (g_M * g_N)(x) \leq C g_L(x)\]
for all $x \in  G$.
\end{lemma}

{\bf Proof} For (a), we use the triangle inequality for $G$:
for some $C > 0$,
\[ |y^{-1}x| \leq C(|y|  +|x|) \]
for all $x,y \in G$. From this we find at once that
\[  (1+|y^{-1}x|)^N \leq C^N(1+|y|)^N(1+|x|)^N, \]
and (a) now follows.\\

For  (b), we similarly observe that
\[ (1+|x|)^L \leq C^L(1+|y|)^L(1+|y^{-1}x|)^L.  \]
Accordingly
\[ (1+|x|)^L (g_M * g_N)(x) = \int (1+|x|)^L  g_M(y)g_N(y^{-1}x)dy \leq
C^L \int g_{M-L}(y) g_{N-L}(y^{-1}x)dy \leq
C^L\|g_{M-L}\|_2 \|g_{N-L}\|_2 \]
which is finite, since $M-L, N-L >  Q/2$.  This completes the proof.\\

Note that Lemma \ref{gmlem} (a) implies that for any measurable subset $E
\subseteq B_0$ of positive measure, we have that
\[ g_N(x) \leq  \frac{C}{m(E)}\int_E g_N(y^{-1}x)dy. \]

for all $x \in G$.  Such facts  will be used without further comment in the
proof which follows.

{\bf Proof of Theorem \ref{framebdd}} We first prove (a). Since we shall be
using Theorem \ref{tof1q}
in our proof of (c), we shall  actually need a stronger conclusion than (a).

We shall in fact show that:\\
\ \\
(*)  For all normalized bump functions $f$ and all $R > 0$ and $u \in G$,
there  exists $C > 0$
such that the series defining $S_{\phi,\psi,b}f^{R,u}$
converges absolutely, uniformly on $G$, and $\|S_{\phi,\psi,b}f^{R,u}\|_{\infty}
\leq C/b^Q$.\\
\ \\
We begin by noting that there exists $C > 0$ such that for any $f, R, u$ 
as above we have:
\[ |<f^{R,u}, \phi_{j,b\gamma}>| \leq  \|f^{R,u}\|_1\|\phi_{j,b\gamma}\|_{\infty}
 \leq C R^Q/a^{jQ/2}. \]

We let

\[ C_{j,R} = \sup |<f^{R,u}, \phi_{j,b\gamma}>|, \]
the sup being taken over all normalized bump functions $f$, all $R > 0$ and all
$u \in G$, and all $\gamma \in \Gamma$.  
Thus

\begin{equation}
\label{frufst}
C_{j,R} \leq C R^Q/a^{jQ/2}.
\end{equation}

Note also that if $R \geq a^j$, then

\begin{equation}
\label{fruscnd}
C_{j,R} \leq C a^{j(Q/2+1)}/R. 
\end{equation}

Indeed, say $f, R, u, \gamma$ are as above.
Since $\int \phi = 0$, putting $v = b\gamma$ we have by Lemma \ref{nbfmvt} that
\begin{eqnarray*}
|<f^{R,u}, \phi_{j,b\gamma}>|   & = &
a^{-jQ/2}|\int_G f^{R,u}(y) \phi( a^{-j}[(a^j v^{-1})y]])  dy|\\
&=& a^{-jQ/2}|\int_G f^{R,u}((a^j v)y) \phi( a^{-j}y)  dy|\\
&=& a^{-jQ/2}|\int_G [f^{R,u}((a^j v)y) -   f^{R,u}(a^jv)] \phi( a^{-j}y)
dy|\\
&\leq& C a^{-jQ/2}
\int_G (\sum_{k=1}^n \frac{|y_k|}{R^{a_k}}) |\phi( a^{-j}y)| dy\\
&=& C  a^{jQ/2} \sum_{k=1}^n [\frac{a^j}{R}]^{a_k} \int_G |y_k \phi(y)| dy\\
&\leq& C  a^{j(Q/2+1)}/R,
\end{eqnarray*}

since we are here
assuming that $a^j/R \leq 1$.

Now select any $N > Q+1$, and note that  $|\psi| \leq Cg_N$ for some $C$.
Fixing $j \in \ZZ$, we now see that

\begin{eqnarray*}
\sum_{\gamma \in \Gamma}|<f^{R,u}, \phi_{j,b\gamma}> \psi_{j,b\gamma}(x) |  &\leq& C_{j,R} C \sum_{\gamma \in \Gamma} D_{a^j}T_{b\gamma}  g_N(x)\\
&=& C a^{-jQ/2}C_{j,R} \sum_{\gamma \in \Gamma}  g_N([b\gamma]^{-1}[a^{-j}x])\\
&\leq& C\frac{a^{-jQ/2} C_{j,R}}{b^Q}\sum_{\gamma \in \Gamma}\int_{b\mathcal{R}} g_N(y^{-1}[b\gamma]^{-1}[a^{-j}x])dy \\
&=&C\frac{a^{-jQ/2} C_{j,R}}{b^Q}\int_G g_N(z)dz\\
&=&C\frac{a^{-jQ/2} C_{j,R}}{b^Q}.
\end{eqnarray*}

Given $R > 0$, we now select $j_0 \in  \ZZ$ with $a^{j_0} \leq R \leq
a^{j_0+1}$.  Recalling (\ref{frufst}) and (\ref{fruscnd}), we now obtain

\begin{eqnarray*}
\sum_{\gamma \in \Gamma, j \in  \ZZ}
|<f^{R,u}, \phi_{j,b\gamma}> \psi_{j,b\gamma}(x)| 
&\leq&
 \frac{C}{b^Q}[\sum_{j \leq j_0}C_{j,R}a^{-jQ/2}
+  \sum_{j > j_0}C_{j,R}a^{-jQ/2}]\\
&\leq&
\frac{C}{b^Q}[\sum_{j \leq  j_0}\frac{a^j}{R}
+ \sum_{j >  j_0}R^Qa^{-jQ}]\\
&\leq&
\frac{C}{b^Q}[\sum_{j \leq  j_0}a^{j-j_0}
+ \sum_{j >  j_0}a^{(j_0-j+1)Q}]\\
&\leq&
\frac{C}{b^Q},
\end{eqnarray*}

proving  (*) and hence (a). \\

We next prove (b).  Again, we shall prove a stronger conclusion, which
will be needed in our proof of (c). \\

For $x,y \in G$, $x \neq y$, we wish to define
\begin{equation}\label{Kphpsdf}
K_{\phi,\psi,b}(x,y) = \sum_{\gamma \in \Gamma, j \in \ZZ}
\psi_{j,b\gamma}(x)\overline{\phi_{j,b\gamma}}(y);
\end{equation}
we will soon show that the sum converges absolutely.  The reason
for this definition is that formally
\[ [S_{\phi,\psi,b} f](x) = \int_G K_{\phi,\psi,b}(x,y)f(y) dy; \]
we will soon show that this is true if $f \in C_c^{1}$, for $x$
outside the support of $f$.  These facts are immediate consequences
of the following assertion (with $J = 0$, $\psi = \Phi$ and $\phi =
\overline{\Psi}$:)

\ \\

(**) 
Say $\Phi, \Psi \in \mathcal{S}(G)$, and $J > 0$.  Then for some $C >0$,

\begin{equation}\label{minQJ}
\sum_{\gamma \in \Gamma, j \in \ZZ} a^{-jJ}|\Phi_{j,b\gamma}(x) \Psi_{j,b\gamma}(y)| 
\leq \frac{C}{b^Q}|y^{-1}x|^{-Q-J}
\end{equation}
for all $x,y \in G$, $x \neq y$.  Moreover the series converges uniformly
on compact subsets of $(G \setminus \{0\}) \times (G \setminus \{0\})$.\\

To prove (\ref{minQJ}), define
\begin{eqnarray*}
K_I(x,y)  &=&
\sum_{\gamma \in \Gamma, j \in \ZZ} a^{-jJ}|y^{-1}x|^{Q+J}
|\Phi_{j,b\gamma}(x)
\Psi_{j,b\gamma}(y)| \\
&=&
\sum_{\gamma \in \Gamma, j \in \ZZ} a^{-j(Q+J)}|y^{-1}x|^{Q+J} 
|\Phi((b\gamma)^{-1}[a^{-j}x])
\Psi((b\gamma)^{-1}[a^{-j}y])|.
\end{eqnarray*}

We need to show that $K_I$ is bounded for $x \neq y$.  Observe that for
any $x,y \in G$ we have that $K_I(ax,ay) = K_I(x,y)$.  Therefore we need
only consider those $x,y$ with $|y^{-1}x| \in [1,a]$.  
Choose $L > Q+J$ and $N > Q/2 + L$.
For some $C_0$, $\Phi, \Psi \leq C_0g_N$.  Thus, for fixed $j$,

\begin{eqnarray*}
\sum_{\gamma \in \Gamma}
|\Phi((b\gamma)^{-1}[a^{-j}x]) \Psi((b\gamma)^{-1}[a^{-j}y])|
&\leq&
C\sum_{\gamma \in \Gamma}
g_N((b\gamma)^{-1}[a^{-j}x])
g_N((b\gamma)^{-1}[a^{-j}y])\\
&\leq&
\frac{C}{b^Q}\sum_{\gamma \in \Gamma}\int_{b\mathcal{R}}
g_N(z^{-1}(b\gamma)^{-1}[a^{-j}x])
g_N(z^{-1}(b\gamma)^{-1}[a^{-j}y])\\
&=&
\frac{C}{b^Q}\int_G g_N(w^{-1}[a^{-j}x]) g_N(w^{-1}[a^{-j}y])dw\\
&=&
\frac{C}{b^Q} (g_N*g_N)(a^{-j}[y^{-1}x])\\
&\leq&
\frac{C}{b^Q} g_L(a^{-j}[y^{-1}x]).
\end{eqnarray*}

Consequently, for $|y^{-1}x| \in [1,a]$, we have that

\begin{eqnarray*}
|K_I(x,y)| &\leq&
\frac{Ca^{Q+J}}{b^Q}\sum_{j \in \ZZ} a^{-j(Q+J)}g_L(a^{-j}[y^{-1}x])\\
&\leq&
\frac{C}{b^Q}[\sum_{j \geq 0} a^{-j(Q+J)} +
\sum_{j < 0} a^{-j(Q+J)}|a^{-j}[y^{-1}x]|^{-L}]\\
&=&
\frac{C}{b^Q}[\sum_{j \geq 0 } a^{-j(Q+J)}
+ \sum_{j < 0} a^{j(L-Q-J)}]\\
&\leq&
\frac{C}{b^Q},
\end{eqnarray*}
proving (\ref{minQJ}).  The uniform convergence asserted in (**) follows
from an examination of the proof of (\ref{minQJ}).  This proves (**).

(**) now implies at once that the series in (\ref{Kphpsdf}) converges absolutely
for $x \neq y$, and we define its sum to be $K_{\phi,\psi,b}(x,y)$.

We can now easily prove (b).  Actually, in order to also later prove (c), we will
note the following stronger statement:\\
\ \\
(***) $K_{\phi,\psi,b}$ is smooth away from the diagonal; moreover for all
multiindices $\alpha, \beta$ there exists $C_{\alpha,\beta}>0$ such that for all
$x,y \in G$ with $x \neq y$ and for all $0 < b < 1$ we have

\begin{equation} \label{xxxyk1}
|X_x^{\alpha} X_y^{\beta}K_{\phi,\psi,b}(x,y)| \leq
\frac{C_{\alpha,\beta}}{b^Q}|y^{-1}x|^{-Q-|\alpha|-|\beta|}.
\end{equation}

(***) follows at once from (**).  Indeed, we claim that, if $x \neq y$, then
\begin{equation}\label{justpull}
X_x^{\alpha} X_y^{\beta}K_I(x,y) =
\sum_{\gamma \in \Gamma, j \in \ZZ} a^{-j(|\alpha|+|\beta|)}
(X^{\alpha}\psi)_{j,b\gamma}(x) (\overline{X^{\beta}\phi})_{j,b\gamma}(y).
\end{equation}
\ \\
To see this, note that by (**), the series in (\ref{justpull}) converges 
uniformly on compact subsets of
\newline
 $(G \setminus \{0\}) \times (G \setminus \{0\})$.
On such a compact set, any (usual) differential monomial 
$\partial_x^{\rho} \partial_y^{\tau}$ is a linear combination, with polynomial
coefficients, of the $X_x^{\alpha} X_y^{\beta}$.  Thus the series 
in (\ref{Kphpsdf}) converges in the topology of 
\newline
$C^{\infty}((G \setminus \{0\}) \times (G \setminus \{0\}))$
This implies that $K$ is smooth off the diagonal, and also that, when
we differentiate $K$, we can bring derivatives past the summation sign.  
This proves (***).\\
  
We now show that
(***) implies (b).
In fact it implies the following stronger statement, which we shall
also need in the proof of (c):\\
\ \\
(****)
For all normalized bump functions $f$ and all $R > 0$ and $u \in G$,
there  exists $C > 0$ such that
\[ \|S_{\phi,\psi,b}f^{R,u}\|_2 \leq \frac{C}{b^Q} R^{Q/2}. \]
\ \\
To see this, 
we observe that if $|x^{-1}u| \geq 2R$, then

\begin{eqnarray}
|[S_{\phi,\psi,b}f^{R,u}](x)| &\leq& \int_G |K_{\phi,\psi,b}(x,y)f^{R,u}(y)|dy \nonumber \\
&\leq& \frac{C}{b^Q}\int_{|u^{-1}y| \leq R} |x^{-1}y|^{-Q}dy \nonumber \\
&\leq& \frac{CR^Q}{b^Q}\max_{|u^{-1}y| \leq R} |x^{-1}y|^{-Q} \nonumber \\
&\leq& \frac{CR^Q}{b^Q} |x^{-1}u|^{-Q}. \label{nbfnegQ}
\end{eqnarray}

(The last inequality follows from Proposition \ref{triget}.)
Finally, if $g = S_{\phi,\psi,b}f^{R,u}$, recalling (*), we have that
\begin{eqnarray*}
\|g\|_2^2 &=& \int_{|x^{-1}u|<2R} |g(x)|^2 dx
+ \int_{|x^{-1}u|>2R}|g(x)|^2 dx\\
&\leq&
\frac{C}{b^{2Q}}[(2R)^Q + R^{2Q}\int_{|x^{-1}u|>2R}|x^{-1}u|^{-2Q}dx]\\
&\leq&
\frac{C}{b^{2Q}}[(2R)^Q + R^{2Q}\int_{|y|>2R}|y|^{-2Q}dy]\\
&\leq&
\frac{C}{b^{2Q}}[(2R)^Q + R^{2Q}(2R)^{-Q}]\\
&=&
\frac{CR^Q}{b^{2Q}}.
\end{eqnarray*}

This proves (****), and hence (b).

We now claim that (c) follows directly from (***), (****) and Theorem \ref{tof1q}. 
In order to apply Theorem \ref{tof1q}, we make the following two
additional observations.
\begin{enumerate}
\item
{\em The formal adjoint of $S_{\phi,\psi,b}$ 
is $S_{\psi,\phi,b}$.}  What we are claiming is that for 
all $f, g \in C_c^{1}$, we have
\begin{equation}\label{eq:49}
<S_{\phi,\psi,b}f,g> = <f,S_{\psi,\phi,b}g>.
\end{equation}
Indeed, by (*), the
left side of (\ref{eq:49}) clearly equals
\[ \sum_{\gamma\in\Gamma, j\in\ZZ}
<f, \phi_{j,b\gamma}> <\psi_{j,b\gamma},g>. \]
Evidently this equals the right side of (\ref{eq:49}), as claimed.
Note, for later purposes, 
that this observation also proves (d) if $f,g \in C_c^{1}$.
\item
{\em $S_{\phi,\psi,b}(1) = S_{\psi,\phi,b}(1) = 0.$} 
We need to show that whenever $f \in C_c^{1}(G)$
and $\int f = 0$, we have that
\[ \int S_{\phi,\psi,b}f = 0 \]
(and similarly $\int S_{\psi,\phi,b}f = 0$.)  To see this, for any finite subset
$\mathcal{F}$ of $\ZZ \times \Gamma$, 
define the operator 
\[ S_{\phi,\psi,b}^{\mathcal{F}}:\; f\rightarrow  
\sum_{(j,\gamma) \in {\cal F}} 
<f, \phi_{j,b\gamma}> \psi_{j,b\gamma}. \]
We regard this as an operator on $C_c^{1}(G)$, and it maps 
this space into $C^{\infty}(G)$, since it has smooth kernel

\[ K_{\phi,\psi,b}^{\mathcal{F}}(x,y) = 
\sum_{(j,\gamma) \in \mathcal{F}}
\psi_{j,b\gamma}(x) \overline{\phi}_{j,b\gamma}(y). \]
For any integer $N > 0$
we also let $S_{\phi,\psi,b}^{N} = S_{\phi,\psi,b}^{\mathcal{F}_N}$
and $K_{\phi,\psi,b}^{N} = K_{\phi,\psi,b}^{\mathcal{F}_N}$, where 
$\mathcal{F}_N = \{(j,\gamma): |j| \leq N, |\gamma| \leq N\}$.
Since $\psi \in \mathcal{S}$ has integral zero, it is evident that for all 
$f \in C_c^{1}(G)$, we have 
\[ \int S_{\phi,\psi,b}^{\mathcal{F}} f = 0 \]
for all $\mathcal{F}$.
Fix $f$ with $\int f = 0$, and fix $b > 0$; we need to deduce 
that $\int S_{\phi,\psi,b} = 0$.  This will follow at once from the dominated covergence
theorem
if we can show:\\

\ \\
(i) $S_{\phi,\psi,b}^N f \rightarrow S_{\phi,\psi,b}f$ pointwise as $N \rightarrow \infty$; and\\

(ii) For some $C > 0$, $|S_{\phi,\psi,b}^{\mathcal{F}} f| \leq C g_{Q+1}$ for all 
$\mathcal{F}$.  \\
\ \\
(Here $g_{Q+1}$ is as in Lemma \ref{gmlem}.)  Since (i) follows at once from the
absolute convergence proved in (*), we need only establish (ii).  (*) similarly
shows that, for some $C > 0$, $|[S_{\phi,\psi,b}^\mathcal{F} f](x)| \leq C$ for all $\mathcal{F}$ and 
all $x \in G$.  Suppose then that the support of $f$ is contained in $\{x: |x| < R\}$;
we need only show that for some $C, A$ (independent of $\mathcal{F}$), 
\begin{equation} \label{zermom}
|[S_{\phi,\psi,b}^\mathcal{F} f](x)| \leq C |x|^{-Q-1} 
\end{equation}
whenever $|x| > AR$.
  But by (**), for any multiindices $\alpha,\beta$, there is a
$C_{\alpha,\beta} > 0$ (independent of $\mathcal{F}$) such that
for all $x, y \in G$ with $x \neq y$,
\[ |X_x^{\alpha}X_y^{\beta}K_{\phi,\psi,b}^\mathcal{F}(x,y)| \leq
C_{\alpha,\beta} |y^{-1}x|^{-Q-|\alpha|-|\beta|}. \]
By Proposition \ref{czmvt} (and its proof), the $K_{\phi,\psi,b}^\mathcal{F}$
satisfy the Calderon-Zygmund inequalities (\ref{cz1}), (\ref{cz2}),
and (\ref{cz3}) with constants $c, C$ independent of $\mathcal{F}$.  
By the triangle inequality, there is a number $A > 0$ such that
whenever $|x| > AR$ and $|y| < R$, we have $|y^{-1}x| > cR$.  
Thus, if $|x| > AR$, we have that
\[ |[S_{\phi,\psi,b}^\mathcal{F} f](x)| = 
|\int_{|y|<R}[K_{\phi,\psi,b}^\mathcal{F}(x,y) - K_{\phi,\psi,b}^\mathcal{F}(x,0)]f(y)dy|
\leq C|x|^{-Q-1}\int_{|y| < R} |y||f(y)|dy \]
as claimed.
\end{enumerate} 

These observations now prove (c) at once. 

We next prove (d).  We fix $b > 0$.
In observation \#1 above, we have already seen that (d) holds for $f,g \in C_c^{1}$.
To prove it in general, we let $S_{\phi,\psi,b}^\mathcal{F}$, $K_{\phi,\psi,b}^\mathcal{F}$,
$S_{\phi,\psi,b}^N$, and $K_{\phi,\psi,b}^N$ be as 
in observation \#2.  We observe that, for any $f,g \in L^2$, we have
\begin{equation}\label{eq:wkN}
<S_{\phi,\psi,b}^N f,g>= 
\sum_{\gamma\in\Gamma, j\in\ZZ,|\gamma|\leq N, |j| \leq N}<f,\phi_{j,b\gamma}>
<\psi_{j,b\gamma},g>,
\end{equation}
and we claim that
\begin{equation}\label{eq:50}
<S_{\phi,\psi,b}^N f,g> \rightarrow <S_{\phi,\psi,b} f,g>.  
\end{equation}
Since (d) holds for $f,g \in C_c^{1}$, which is dense in $L^2$,
it is enough to show that the norms $\|S_{\phi,\psi,b}^\mathcal{F} \|$ are uniformly bounded
in $\mathcal{F}$.  But this follows from Theorem \ref{tof1q}, together with a repetition of the proofs of
(*), (***) and (****) with $S_{\phi,\psi,b}^\mathcal{F}$, $K_{\phi,\psi,b}^\mathcal{F}$ in place of
$S_{\phi,\psi,b}$, $K_{\phi,\psi,b}$; here one must note that all bounds
are independent of $\mathcal{F}$.  

Now in (\ref{eq:50}), take the special case $\psi = \phi$ and $g = f$.  
All the terms in the series in 
(\ref{eq:wkN}) are then nonnegative, so the series in (\ref{eq:wk}) converges absolutely

to the left side of that equation -- in that special case.  But in the general case,
Cauchy-Schwarz as applied to the series in (\ref{eq:50}) now shows that this series
always converges absolutely.  Moreover, by (\ref{eq:50}), this series converges to the
left side of (\ref{eq:wk}).  This proves (d).

(e) follows from an examination of the proofs
of (a).  (f) follows from an examination of the proofs of (a), (b) and (c).  
(In particular, note for later purposes that the constants $C_{\alpha,\beta}$ in 
(\ref{xxxyk1}) may be taken independent of $\phi, \psi \in \BB_0$ with $\int \psi
= \int \phi = 0$.)  This completes the proof of Theorem \ref{framebdd}.
\\

The main idea in our proof of Theorem \ref{T:3} (a) is to show that, for $b$ sufficiently 
small, $Vb^Q S_{\psi,\psi,b}$ is ``well approximated'' by the operator 
$R_{\psi} = \sum_{j \in \ZZ} R_j$,
where if $f \in L^2(G)$ we put
\[ R_j f = f * \tilde{\psi}_{a^j} * \psi_{a^j},  \]
then to use the spectral theorem to show that $R_{\psi}$ is bounded below if $\psi = g(L)\delta$
and $g$ satisfies Daubechies' criterion.  We begin to make these ideas rigorous, by noting 
the following proposition.  (In this proposition, $\mathcal{R}$ is, once again,
a bounded measurable subset of $G$, of positive measure, such that every
$g \in G$ may be written uniquely
in the form $g = x\gamma$ with $x \in \mathcal{R}$ and $\gamma \in \Gamma$.)

\begin{proposition}\label{P:1}
Suppose $\psi \in \mathcal{S}(G)$ and $\int \psi = 0$.  Suppose $f \in C_c^{1}(G)$.
Then 
\begin{equation}\label{eq:51}
\sum_{j \in \ZZ} (R_j f)(x) = \int_{b\mathcal{R}} [S_{T_z\psi,T_z\psi,b}f](x)dz
\end{equation}

where the series on the left side converges absolutely, uniformly for $x$ on G. 
Consequently $\sum R_j f$
converges to an $L^2$ function, and the map $R_{\psi}: C_c^{1} \rightarrow L^2$
given by $R_{\psi}f = \sum R_j f$ extends to a bounded positive operator on $L^2$.\\
\end{proposition}
{\bf Proof}  Fix $j$ for now and put $\eta = \psi_{a^j}$.  Then 
\[(R_j f)(x) = \int_G f(y)[\tilde{\eta}*\eta](y^{-1}x) dy. \]
But
\begin{eqnarray*}
[\tilde{\eta}*\eta](y^{-1}x) dy &=& 
\int_G \tilde{\eta}(y^{-1}z) \eta(z^{-1}x) dz\\
&=&
\int_G \overline{\eta}(z^{-1}y) \eta(z^{-1}x) dz\\
&=&
a^{-2jQ}\int_G \overline{\psi}(a^{-j}[z^{-1}y]) \psi(a^{-j}[z^{-1}x]) dz\\
&=&
a^{-jQ}\int_G \overline{\psi}(z^{-1}[a^{-j}y]) \psi(z^{-1}[a^{-j}x]) dz\\
&=&
a^{-jQ}\sum_{\gamma \in \Gamma} 
\int_{b\mathcal{R}} 
\overline{\psi}(z^{-1}[b\gamma]^{-1}[a^{-j}y]) \psi(z^{-1}[b\gamma]^{-1}[a^{-j}x]) dz\\
&=&
\sum_{\gamma \in \Gamma} \int_{b\mathcal{R}} 
\overline{(T_z \psi)}_{j,b\gamma}(y) (T_z \psi)_{j,b\gamma}(x)dz.\\
\end{eqnarray*}
(In the fourth line, we have made the change of variables $z \rightarrow a^jz$.)
By Theorem \ref{framebdd} (e), the series
\[ \sum_{j \in \ZZ, \gamma \in \Gamma} <f,(T_z \psi)_{j,b\gamma}>(T_z \psi)_{j,b\gamma}(x) \]
converges absolutely, uniformly for $x \in G$ and $z \in b\mathcal{R}$.  This would therefore
also surely be true if we fixed a $j$ and summed only over $\gamma$.  Thus
\[ (R_j f)(x) = \int_{b\mathcal{R}} \sum_{\gamma \in \Gamma} 
<f,(T_z \psi)_{j,b\gamma}>(T_z \psi)_{j,b\gamma}(x)dz \]
and finally, summing over $j$, we find (\ref{eq:51}) as well, the sum on the left side
converging absolutely, uniformly for $x \in G$.
The remaining conclusions of the proposition now follow
at once from Theorem \ref{framebdd} (f) and Minkowski's inequality.  (Note: to show
that $R_{\psi}$ is positive, it is enough to show that $<R_{\psi}f,f> \geq 0$ for all $f \in
C_c^{1}$, and for this it is enough to show that $<R_j f,f> \geq 0$ for all
$f \in C_c^{1}$ and all $j$.  But this is clear, since for such $f$,
$<R_j f,f> = \|f*\psi_{a^j}\|_2^2$.)
 This completes the proof.\\

We can now reach an understanding of why
$Vb^Q S_{\psi,\psi,b}$ is well approximated by $R_{\psi}$ for $b$ small.

\begin{theorem}\label{genapp}
Suppose $\psi \in \mathcal{S}(G)$, and $\int_G \psi = 0$.  Let $V$ be the 
measure of $\mathcal{R}$, and let $R_{\psi}$ be as in Proposition \ref{P:1}.
For $1 \leq l \leq n$, let $\psi_l = Y_l \psi$ (so that, by Proposition
\ref{zeroint}, $\int_G \psi_l = 0$
for all $l$).  Then:\\
(a) If $f \in C_c^{1}(G)$, 
\begin{equation}\label{eq:53}
[\frac{1}{Vb^Q} R_{\psi} f - S_{\psi,\psi,b} f  ](x) 
=\frac{1}{Vb^Q} \sum_{l=1}^n \int_{b \mathcal{R}} \int_0^1
z_l ([S_{T_{[t]z}Y_l \psi, T_{[t]z}\psi,b}f] + 
[S_{T_{[t]z}\psi, T_{[t]z}Y_l \psi,b}f])(x) dt dz.
\end{equation}
(b) There exists $C > 0$ such that for all $0 < b < 1$, the norm on $L^2(G)$
\begin{equation}\label{eq:54}
\| \frac{1}{Vb^Q} R_{\psi} - S_{\psi,\psi,b}\| \leq \frac{C}{b^{Q-1}}. 
\end{equation}
(c) If $R_{\psi} \geq AI$ for some $A > 0$, then there exists $b_0 > 0$ such
that $\{\psi_{j,b\gamma}\}$ is a frame whenever $0 < b < b_0$. 
More precisely, choose $B > 0$ such that $R_{\psi} \leq BI$ (of course we
can choose $B = \|R_{\psi}\|$).  Then, for $0 < b < b_0$, we can choose
$A_b, B_b > 0$ such that 
\begin{equation}\label{eq:55}
A_b\|f\|_2^2  
\leq \sum_{j \in \ZZ, \gamma \in \Gamma}|<f,\psi_{j,b\gamma}>|^2
\leq B_b\|f\|_2^2 
\end{equation}
for all $f \in L^2$, and such that
\begin{equation}\label{eq:56}
\lim_{b \rightarrow 0^+} \frac{B_b}{A_b} = \frac{B}{A}.
\end{equation}
\end{theorem}
{\bf Proof}  Of course, the measure of $b\mathcal{R}$ is $Vb^Q$.
Using Theorem \ref{framebdd} (e), together with Proposition \ref{P:1}, we see that
\begin{eqnarray}
[\frac{1}{Vb^Q} R_{\psi} f - S_{\psi,\psi,b} f](x) &=&
\frac{1}{Vb^Q}\int_{b\mathcal{R}}([S_{T_z\psi,T_z\psi,b}f](x)-[S_{\psi,\psi,b}f](x))dz\\
&=&
\frac{1}{Vb^Q}\sum_{j \in \ZZ, \gamma \in \Gamma}\int_{b\mathcal{R}} 
[<f,(T_z \psi)_{j,b\gamma}>(T_z \psi)_{j,b\gamma}(x) 
- <f,\psi_{j,b\gamma}>\psi_{j,b\gamma}(x)]dx \nonumber \\
\ & \ & \ \ \label{eq:dif}
\end{eqnarray}
However, fixing $j, \gamma$, we have that
\begin{equation}\label{eq:54a}
<f,(T_z \psi)_{j,b\gamma}>(T_z \psi)_{j,b\gamma}(x)
- <f,\psi_{j,b\gamma}>\psi_{j,b\gamma}(x) =
\int_G f(y) K(x,y) dy
\end{equation}
where 
\[ K(x,y) = \overline{(T_z \psi)}_{j,b\gamma}(y)(T_z \psi)_{j,b\gamma}(x)
- \overline{\psi}_{j,b\gamma}(y)\psi_{j,b\gamma}(x). \]
Fix $x,y$ as well, and, for $w \in G$, let 
\[ F(w) = \overline{(T_{w^{-1}} \psi)}_{j,b\gamma}(y)(T_{w^{-1}} \psi)_{j,b\gamma}(x). \]
Explicitly
\[ F(w) =  a^{-jQ}\overline{\psi}(w(b\gamma)^{-1}(a^{-j}y))\psi(w(b\gamma)^{-1}(a^{-j}x)).\]
Since each $Y_l$ is right-invariant, note that
\begin{eqnarray*}
(Y_1 F)(w) &=&
 a^{-jQ}[\overline{(Y_l\psi)}(w(b\gamma)^{-1}(a^{-j}y))\psi(w(b\gamma)^{-1}(a^{-j}x))+
\overline{\psi}(w(b\gamma)^{-1}(a^{-j}y))(Y_l\psi)(w(b\gamma)^{-1}(a^{-j}x))]\\ 
&=&
\overline{(T_{w^{-1}} [Y_l\psi])}_{j,b\gamma}(y)(T_{w^{-1}} \psi)_{j,b\gamma}(x) +  \overline{(T_{w^{-1}} \psi)}_{j,b\gamma}(y)(T_{w^{-1}} [Y_l\psi])_{j,b\gamma}(x).
\end{eqnarray*} 
Then, by Lemma \ref{ftc}, and the facts that $z^{-1} = -z$ and
that the $Y_l$ are right-invariant, we have
\begin{eqnarray*}
K(x,y) &=& F(z^{-1}) - F(0) \\ \nonumber
&=& -\sum_{l=1}^n z_l \int_0^1(Y_lF)({([t]z)}^{-1})dt \\ 
&=& -\sum_{l=1}^n z_l \int_0^1 
[\overline{(T_{[t]z} [Y_l\psi])}_{j,b\gamma}(y)(T_{[t]z} \psi)_{j,b\gamma}(x) + 
\overline{(T_{[t]z} \psi)}_{j,b\gamma}(y)(T_{[t]z} [Y_l\psi])_{j,b\gamma}(x)]dt. 
\end{eqnarray*}
Since $f \in C_c^{1}$,

\[ \int_G f(y)K(x,y)dy =  -\sum_{l=1}^n z_l \int_0^1 
[<f,{(T_{[t]z} [Y_l\psi])}_{j,b\gamma}>(T_{[t]z} \psi)_{j,b\gamma}(x) + 
<f,{(T_{[t]z} \psi)}_{j,b\gamma}>(T_{[t]z} [Y_l\psi])_{j,b\gamma}(x)]dt. \]

Part (a) of the theorem now follows at once from this, (\ref{eq:dif}),
(\ref{eq:54a}), and Theorem \ref{framebdd} (e).

For (b), choose a number $M > 0$ such that $|z_l| \leq M$ whenever $z \in \mathcal{R}$
and $1 \leq l \leq n$. Then, surely, $|z_l| \leq b^{a_l}M \leq bM$ whenever
$z \in b\mathcal{R}$ and $1 \leq l \leq n$.  Accordingly, by
Theorem \ref{framebdd} (f) and Minkowski's inequality,
there exist $C_1, C > 0$ such that for all $f \in C_c^{1}$ and all $0 < b < 1$, we
have
\begin{eqnarray*}
\|\frac{1}{Vb^Q}R_{\psi} f - S_{\psi,\psi,b}] f  ]\|_2 
&\leq &
\frac{bM}{Vb^Q} \sum_{l=1}^n \int_{b \mathcal{R}} \int_0^1
(\|[S_{T_{[t]z}Y_l \psi, T_{[t]z}\psi,b}f\|_2 + 
\|S_{T_{[t]z}\psi, T_{[t]z}Y_l \psi,b}f\|_2) dt dz.\\
&\leq& \frac{2nbM}{Vb^Q}m(b\mathcal{R}) \frac{C_1}{b^Q}\|f\|_2\\
&=& \frac{C}{b^{Q-1}}\|f\|_2.
\end{eqnarray*}
Since $C_c^{1}$ is dense in $L^2$, (b) now follows.

Finally for (c), take $C$ as in (b).  Then for any $f \in L^2$,
\[ <S_{\psi,\psi,b}f,f> = \frac{1}{Vb^Q}<R_{\psi}f,f> + 
<[S_{\psi,\psi,b} - \frac{1}{Vb^Q}R_{\psi}]f,f> \]
so that, for all $f \in L^2$,
\[ A_b\|f\|_2^2\ \  \leq\ \  <S_{\psi,\psi,b}f,f>\ \  \leq \ \  B_b\|f\|_2^2,  \]
where
\begin{equation}\label{eq:57}
A_b = \frac{A-CVb}{Vb^Q}
\end{equation}
and 
\begin{equation}\label{eq:58}
B_b = \frac{B+CVb}{Vb^Q}.
\end{equation}

 Put $b_0 = \min(A/CV,1)$.  If $0 < b < b_0$, then $A_b > 0$, and, moreover,
by Theorem \ref{framebdd} (d), (\ref{eq:55}) holds for all $f \in L^2$.  
Finally (\ref{eq:56}) is immediate from (\ref{eq:57}) and (\ref{eq:58}).
This proves (c) and completes the proof of the theorem.  \\

Accordingly, the search for frames reduces to the question of finding $\psi$
with $R_{\psi} \geq AI$ for some $A > 0$.  \\

If $G = \RR^n$ with the usual addition, then if $f \in L^2$,
\[  \widehat{(R_j f)}(\xi) = |\hat{\psi}(a^j \xi)|^2 \hat{f}(\xi)  \]
Define
\[  m_{\hat{\psi},a}(\xi) = \sum_{j \in \ZZ} |\hat{\psi}(a^j \xi)|^2 \]
(Usually $a$ is fixed and understood, and we will just write $m_{\hat{\psi}} =
m_{\hat{\psi},a}$.)  
Since we are assuming that $\int \psi = 0$, surely $\hat{\psi}(0) = 0$,
so that $|\hat{\psi}(\xi)| \leq C|\xi|$ for $|\xi| < 1$. Also, since 
$\hat{\psi} \in \mathcal{S}$, $|\hat{\psi}(\xi)| \leq C/|\xi|$ for $|\xi| \geq 1$.
From these facts it is easy to see that the series defining $m_{\hat{\psi}}(\xi)$ converges
uniformly on any compact subset of $\RR^n$ which excludes the origin.  
We claim that
\[  \widehat{(R_{\psi}f)}(\xi) = m_{\hat{\psi}}(\xi) \hat{f}(\xi)  \]
for all $f \in L^2$.  This is not hard to see, but since we shall need an 
analogue for general $G$, let us present the argument in detail.

First note that
$m_{\hat{\psi}}(a\xi) = m_{\hat{\psi}}(\xi)$ for all $\xi$, so $m_{\hat{\psi}}$ 
is uniformly bounded on $\RR^n$.  Define an operator $Q: L^2 \rightarrow L^2$ by 
$\widehat{(Qf)}(\xi) = m_{\hat{\psi}}(\xi) \hat{f}(\xi)$; we want to show
that $R_{\psi} = Q$.
For $N > 0$, set $Q_N = \sum_{j=-N}^N R_j$, an operator on $L^2$; then 
$\|Q_N\| \leq \|m_{\hat{\psi}}\|_{\infty}$ for all $N$.  If 
\[ V = \{f \in L^2: \hat{f} = 0 \mbox{ a.e. outside some 
compact subset of } \RR^n \setminus \{0\}\}, \]
then $Q_N f \rightarrow Qf$ in $L^2$ for all $f \in V$.  Since $V$ is dense in $L^2$
and the $\|Q_N\|$ are uniformly bounded, we see that $Q_N f \rightarrow Qf$
for all $f \in L^2$.  However, $Q_N f \rightarrow R_{\psi}f$ pointwise on $\RR^n$
if $f \in C_c^{1}$.  Consequently $Qf = R_{\psi}f$ for all $f \in C_c^{1}$,
and hence for all $f \in L^2$, as claimed.

If we now let 
\begin{equation}\label{eq:59}
B = \sup_{\xi \neq 0} m_{\hat{\psi}}(\xi),\ \ \ A = \inf_{\xi \neq 0} m_{\hat{\psi}}(\xi),
\end{equation}
we now see that 
\[ A\|f\|_2^2 \leq (R_{\psi}f,f) \leq B\|f\|_2^2  \]
for all $f \in L^2$.  Theorem \ref{genapp} then tells us in particular that 
$\{\psi_{j,b\gamma}\}$ is a frame for all sufficiently small $b$, {\em provided
that} $A > 0$.  The condition $A > 0$ is called {\em Daubechies' criterion}.  
(In \cite{Daubechies92}, Daubechies shows that $\{\psi_{j,b\gamma}\}$ is a frame
if $n = 1$ and $\Gamma$ is the integer lattice, if Daubechies' criterion holds.
Her methods are very different from those of this paper -- she uses Plancherel and Parseval.)
Since, for all $\xi$, $m_{\hat{\psi}}(a\xi) = m_{\hat{\psi}}(\xi)$, $A$ is the minimum
value of $m_{\hat{\psi}}$ on the compact annulus $\{\xi: 1 \leq |\xi| \leq a\}$.  
Thus Daubechies' criterion is equivalent to the hypothesis that there does not
exist a nonzero $\xi_0 \in \RR^n$ such that $\hat{\psi}(a^j \xi_0) = 0$ for all integers  
$j$.

Now we turn to general stratified Lie groups $G$.\\

{\bf Proof of Theorem \ref{T:3}} We change notation from the statement of Theorem
\ref{T:3}, writing $H$ in place of $f$.  Thus we
restrict attention to $\psi$ of the form $\psi = F(L)\delta
= LH(L)\delta$, where $F(\lambda) = \lambda H(\lambda)$ and $H, F \in \mathcal{S}(\RR^+)$.
In that case, if $f \in L^2$,
\[ R_j f = F_j(L) f \]
where
\[ F_j(\lambda) = |F(a^{2j}\lambda)|^2.  \]
Define
\[  m_{F,a^2}(\lambda) = \sum_{j \in \ZZ} |F(a^{2j} \lambda)|^2 \]
(Usually $a$ is fixed and understood, and we will just write $m_F =
m_{F,a^2}$.)  
As before, the series defining $m_F$ converges
uniformly on any compact subset of $\RR^{+}$ which excludes the origin.  
We claim that
\begin{equation}\label{eq:60}
R_{\psi}f = m_F(L)f 
\end{equation}
for all $f \in L^2$.  

As before, 
$m_F(a^2\lambda) = m_F(\lambda)$ for all $\lambda$, so $m_F$ 
is uniformly bounded on $\RR^{+}$.  Set $Q = m_F(L)$, and put
$\mathcal{H} = L^2(G)$.
For $N > 0$, set $Q_N = \sum_{j=-N}^N R_j$, an operator on $\mathcal{H}$; then,
by the spectral theorem,
$\|Q_N\| \leq \|m_F\|_{\infty}$ for all $N$.  Let
\[ V = \bigcup_{0 < \epsilon < N < \infty} P_{[\epsilon,N]}\mathcal{H}\]
(Recall 
that the $P_{[a,b]}$ are spectral projectors of $L$.)
Then $Q_N f \rightarrow Qf$ in $\mathcal{H}$ for all $f \in V$.  But, since 
$P_{\{0\}} = 0$,
$V$ is dense in $\mathcal{H}$.
Since the $\|Q_N\|$ are uniformly bounded, it follows that $Q_N f \rightarrow Qf$
for all $f \in \mathcal{H}$.  However, $Q_N f \rightarrow R_{\psi}f$ pointwise on $G$
if $f \in C_c^{1}$.  Consequently $Qf = R_{\psi}f$ for all $f \in C_c^{1}$, and hence for all $f \in L^2$, as claimed.

If we now let 
\begin{equation}\label{eq:61}
B = \sup_{\lambda > 0} m_F(\lambda),\ \ \ A = \inf_{\lambda > 0} m_F(\lambda),
\end{equation}
and again note that $P_{\{0\}} = 0$, we see that 
\[ A\|f\|_2^2 \leq (R_{\psi}f,f) \leq B\|f\|_2^2  \]
for all $f \in L^2$.  Theorem \ref{genapp} then tells us in particular that 
$\{\psi_{j,b\gamma}\}$ is a frame for all sufficiently small $b$, {\em provided
that} $A > 0$ -- in other words, if $F$ satisfies Daubechies' criterion (where
of course we use $a^2$ in place of the $a$ we used on $\RR^n$).  

This establishes Theorem \ref{T:3} (a).  (Here Daubechies' criterion is clearly equivalent to 
the nonexistence of a $\lambda_0 > 0$ such that $F(a^{2j}\lambda_0) = 0$ for all
integers $j$.)

We now prove Theorem \ref{T:3} (b).  By (\ref{eq:56}) we need only show that if $a$ is
close enough to $1$, then $A = \inf_{\lambda > 0} m_{F,a^2}(\lambda) > 0$ (so that
the Daubechies condition holds), and that
\begin{equation}
\label{BAone}
\frac{B}{A} = 1 + O(|a-1|^2 \log|a-1|).
\end{equation}

as $a \rightarrow 1$.
We may assume $a > 1$ (otherwise replace $a$ by $1/a$, and note that
$m_{F,a^2} = m_{F,(1/a)^2}$, and \newline
$O(|a-1|^2 \log|a-1|) = O(|1/a-1|^2 \log|1/a-1|)$.)
However, if $a > 1$, then (\ref{BAone}) follows at once from the following elementary lemma.

\begin{lemma}
\label{elem}
Suppose $H$ is a nonzero element of ${\mathcal S}(\RR^+)$ 
and let $F(s) = sH(s)$. Let $I \in (0,\infty)$ be defined by

\begin{equation}
\label{caldc}
I = \int_0^{\infty} |F(t)|^2 \frac{dt}{t} = \int_0^{\infty} |F(ts)|^2 \frac{dt}{t}
\end{equation}

(for any $s > 0$), as in Calder\'on's reproducing formula.  
Suppose $a > 1$.  Then for all $s > 0$,
\begin{equation}
\label{daubtt}
A(a) \leq \sum_{n=-\infty}^{\infty} |F(a^{2n} s)|^2 \leq B(a) < \infty,
\end{equation}
where, as $a \rightarrow 1$,
\begin{equation}
\label{daubest}
A(a) = \frac{I}{2\log a}(1 - O(|a-1|^2\log|a-1|)),\:\:\: B(a)=
\frac{I}{2\log a}(1 + O(|a-1|^2\log|a-1|)).
\end{equation}
\
\end{lemma}
{\bf Proof}
Define a new function $G: \RR \rightarrow \RR$ by 
\begin{equation} 
\label{f1f}
G(u) = |F(e^u)|^2 = |e^u H(e^u)|^2;
\end{equation}
then 
\begin{equation}
\label{f1prps}
G \in \mathcal{S}(\RR), \mbox{ and } |G(u)| \leq Ke^{-2|u|}
\end{equation}
for some constant $K$.\\

If we put $t = e^u$ in Calder\'on's 
identity (\ref{caldc}), and also write $s = e^v$, that identity becomes
the simpler identity

\begin{equation}
\label{cald1}
\int_{-\infty}^{\infty} G(u + v) du = I 
\end{equation}

(independent of $v$).  If we again put $s = e^v$, and now write 
$a^2 = e^c$, we see that the sum we need to estimate has the simpler form

\begin{equation}
\label{dau2}
\sum_{n=-\infty}^{\infty} |F(a^{2n} s)|^2 =
\sum_{n=-\infty}^{\infty} G(nc + v).
\end{equation}

Since the sum on the right side of (\ref{dau2})
is periodic with period $c$, the sum need only be estimated for $0 \leq v \leq c$.
Since we are letting $a \rightarrow 1^+$, we may assume $0 < c = 2 \log a< 1/e$.

We note that $c\sum_{n=-\infty}^{\infty} G(nc + v)$ is a Riemann sum for 
the integral $\int_{-\infty}^{\infty} G(u + v) du = I$. 

To estimate the difference, we recall the midpoint rule:  Say $f$ is $C^2$ in
a neighborhood of $[a,b]$.  Divide $[a,b]$ into $n$ intervals of equal length $\Delta x =
(b-a)/n$ and let $x_k^*$ be the midpoint of the $k$th interval.  Let

\[ E = |\int_a^b f(x)dx - \sum_{k=1}^n f(x_k^*)\Delta x|. \]

Then
\[ E \leq \frac{1}{24}\|f''\|_{\infty}(b-a)(\Delta x)^2. \]

Thus, there is a constant $P > 0$ such that whenever $0 \leq v \leq c < 1/e$,
and whenever $N > 0$ is an integer,
\begin{eqnarray*}
&&|c\sum_{n=-\infty}^{\infty} G(nc + v) - I|\\
& = & |c\sum_{n=-\infty}^{\infty} G(nc + v) - \int_{-\infty}^{\infty} G(u + v)du|\\
& \leq & |c\sum_{n=-N}^{N} G(nc + v) - \int_{-(Nc+c/2)}^{Nc+c/2} G(u + v)du| + 
c\sum_{|n|>N} G(nc + v) + \int_{|u| > Nc + c/2} G(u + v)du\\
& \leq & \frac{1}{24}\|G''\|_{\infty}[(2N+1)c]c^2 + 2 K e^{-2(N+1)c}e^{2v}\frac{c}{1-e^{-2c}}
+  K e^{-(2N+1)c} e^{2v}\\
& \leq & P(Nc^3+e^{-2Nc}). 
\end{eqnarray*}
Note that for $x > e$, $x \log x > e \log e > 1$.  Since we are assuming $1/c > e$,
there is an integer $N$ with $[\log(1/c)]/c < N < [2\log(1/c)]/c$.  Using such an
$N$ we see that
\[ |\sum_{n=-\infty}^{\infty} |F(a^{2n} s)|^2 - \frac{I}{c}| = 
|\sum_{n=-\infty}^{\infty} G(nc + v) - \frac{I}{c}|
\leq P(2c\log(1/c) + c) \leq 3Pc\log(1/c).  \]

Accordingly $\sum_{n=-\infty}^{\infty} |F(a^{2n} s)|^2$ is between
$(I/c)(1 \pm Qc^2|\log c|)$, where $Q = 3P/I$.  Since $c = 2\log a$, 
and since $\log a/(a-1) \rightarrow 1$ as $a \rightarrow 1^+$,
we have completed the proof of Lemma \ref{elem},  and, with it, the proof of Theorem 
\ref{T:3}.\\

{\bf Example } 
Daubechies (\cite{Daubechies92}, especially page 77 and pages 71-72), calculated
for instance that if $a = 2^{1/3}$,\\ $\psi(x) = c(1-x^2)e^{-x^2/2}$ (for $x \in \RR$;
here $c \neq 0$ could be chosen arbitrarily),
$B = \sup_{\xi > 0} m_{\hat{\psi},a}(\xi),\\ A = \inf_{\xi > 0} m_{\hat{\psi},a}(\xi)$, then 
$B/A=1.0000$ to four significant digits.  
\footnote{Actually, Daubechies took a specific value of $c$, but clearly that is 
irrelevant in computing $B/A$.  Also, in her table
on page 77 of \cite{Daubechies92}, her $B/A$ is larger than 
$\sup_{\xi > 0} m_{\hat{\psi},a}(\xi)/\inf_{\xi > 0} m_{\hat{\psi},a}(\xi)$,
 (see her equations (3.3.19) and (3.3.20)),
but that is an even stronger assertion than the one we are making.}

In that case $\psi$ is a multiple of the second derivative of $e^{-x^2/2}$, so
$\hat{\psi}(\xi) = c'\xi^2 e^{-\xi^2/2}$.  Again $c'$ is nonzero and arbitrary;
let us now take it to be $1$.  If we let $F(\lambda) = 
\hat{\psi}(\sqrt{\lambda}) = \lambda e^{-\lambda/2}$ (essentially 
\footnote{If $G=\RR$ we are of course passing from the spectral resolution of
$d/dx$ to that of $d^2/dx^2$.}
making the change of variables $\lambda = \xi^2$), then 
\[ F(L)\delta = L e^{-L/2}\delta = \Psi, \mbox{ say,} \]
and
\[ m_{F,a^2}(\lambda) = m_{\hat{\psi},a}(\sqrt{\lambda}),\]
so that $\sup_{\lambda > 0} m_{F,a^2}(\lambda) =
 \sup_{\xi > 0} m_{\hat{\psi},a}(\xi) = B$, say, and
$\inf_{\lambda > 0} m_{F,a^2}(\lambda) = 
\inf_{\xi > 0} m_{\hat{\psi},a}(\xi) = A$, say.  By the aforementioned
calculation of Daubechies, $B/A = 1.0000$ to four significant digits. 
Thus, by Theorem \ref{genapp}, we can choose $b_0 > 0$ such that
$\{\Psi_{j,b\gamma}\}$ is a frame whenever $0 < b < b_0$, with frame
bounds $A_b, B_b$ and such that, moreover, $B_b/A_b = 1.0000$ to four
significant digits.\\

$\Psi$ is, up to a constant multiple, a natural generalization of 
the Mexican Hat wavelet to $G$.

\section{Frames in Other Banach Spaces}

In this section we discuss the invertibility of $S_{\psi,\psi,b}$
on other Banach spaces, such as $L^p$ or $H^1$.  Let us clarify which
Banach spaces we can allow.
 \begin{definition}
\label{acceptban}
We call a Banach space $\BB$ of measurable functions on $G$ {\em acceptable} if 
$L^2 \cap \BB$ is dense in $\BB$, $\BB \subseteq \mathcal{S}'$ (continuous
inclusion), and if the following condition holds:\\
There exist $C_0, N > 0$, such that for any $A_0 > 0$, we have the following.
Whenever $T: L^2 \rightarrow L^2$ is linear and satisfies:
\ \\
(i) The operator norm of $T$ on $L^2$ is less than or equal to $A_0$;\\ 
(ii) There is a kernel $K(x,y)$, $C^1$ off the diagonal, 
such that if $f \in C_c^1$,
then for $x$ outside the support of $f$, $(Tf)(x) = \int K(x,y)f(y)dy$; and 
whenever $0 \leq \alpha_1 + \ldots + \alpha_n + \beta_1 + \ldots + \beta_n \leq 1$,
and whenever $x, y \in G$ with $x \neq y$, we have
\begin{equation}
\label{xxxyk}
|X_x^{\alpha} X_y^{\beta}K(x,y)| \leq A_0|y^{-1}x|^{-Q-|\alpha|-|\beta|};\ \ \ \ 
{\rm and}
\end{equation}
(iii) $T^*(1) = 0$;
\ \\
then $T|_{L^2 \cap \BB}$ extends to a bounded operator on $\BB$, with norm 
$\|T\| \leq C_0 A_0$.
\end{definition}

Surely (\cite{CoifmanWeiss1}, \cite{CoifmanWeiss2}) $L^p$ ($1 < p < \infty$) 
and $H^1$ are acceptable Banach spaces.  In this section, we shall show:

\begin{theorem} \label{accept}
Suppose $H \in \mathcal{S}(\RR^+)$, $F(\lambda) = \lambda H(\lambda)$, and that $F$ satisfies
Daubechies' criterion (i.e., that 
$\inf_{\lambda > 0} m_{F,a^2}(\lambda) > 0$).  Let $\psi = F(L)\delta$.
Suppose $\BB$ is acceptable, and that $\psi_{j,b\gamma} \in \BB$
for all $j \in \ZZ$ and $0 < b < 1$.  Then:\\
\ \\
(a) For some $b_0 > 0$, $S_{\psi,\psi,b}$ is invertible on $\BB$ whenever $0 < b < b_0$.\\
\ \\
Suppose now that $0 < b < b_0$.\\
\ \\
(b) Suppose that for some dense subspace $\mathcal{D}$ of $\BB$, the series 
$\sum_{j,\gamma} <f,\psi_{j,b\gamma}> \psi_{j,b\gamma}$ converges uncondiionally
to $S_{\psi,\psi,b}f$ in $\BB$ for all $f \in \mathcal{D}$.  Then this series converges
unconditionally to $S_{\psi,\psi,b}f$ for all $f \in \BB$.  Moreover, if we let
$\phi^{j,b\gamma} = S^{-1}_{\psi,\psi,b}\psi_{j,b\gamma}$, then for any $f \in \BB$,
\begin{equation}
\label{uncon1}
f = \sum_{j,\gamma} <f,\psi_{j,b\gamma}> \phi^{j,b\gamma},
\end{equation}
where the series converges unconditionally to $f$ in $\BB$. 
In particular, the set $\{\phi^{j,b\gamma}\}$ is a complete system in $\BB$ (i.e., the
closure of the linear span of this set is all of $\BB$).\\
(c) The hypotheses, and hence the conclusion, of (b) hold if $\BB = L^p$ ($1 < p < \infty$)
or $H^1$.  Here we may take $\mathcal{D} = C_c^{\infty} \cap \BB$.
\end{theorem}
{\bf Proof}
We retain all the notation of the proofs of Theorems \ref{framebdd} and \ref{genapp}.\\

For (a), by Theorem \ref{framebdd} and Defintion \ref{acceptban},
$S_{\psi,\psi,b}|_{L^2 \cap \BB}$ extends to 
a bounded operator on $\BB$.  Also, by Proposition \ref{P:1}, Theorem \ref{framebdd} (f),
the second last sentence of the proof of Theorem \ref{framebdd} (f), and Definition
\ref{acceptban}, we have that $R_{\psi}|_{L^2 \cap \BB}$ extends to 
a bounded operator on $\BB$.  Further, by Theorem \ref{genapp} (a) and Minkowski's
inequality, there exists $C > 0$ such that for all $0 < b < 1$, the norm on $\BB$
\[\| \frac{1}{Vb^Q} R_{\psi} - S_{\psi,\psi,b}\| \leq \frac{C}{b^{Q-1}}. \]
To prove (a), it suffices to show that $R_{\psi}$ is invertible on 
$\BB$.  Indeed, say this were known.  For (a), it is clearly enough to show
that the operator
\[ L_b = Vb^Q R_{\psi}^{-1} S_{\psi,\psi,b}\]
is invertible on $\BB$ for all sufficiently small $b$.  But this is clear, since 
\[ \|I-L_b\| \leq C\|R_{\psi}^{-1}\| b\]
which is less than $1$ if $b$ is sufficiently small.

So it is enough to show that $R_{\psi}$ is invertible on $\BB$.  
By (\ref{eq:60}), $R_{\psi} \equiv m_F(L)$ on $L^2$.  By Daubechies'
criterion, $1/m_{F} = G$, say, is a bounded function on $\RR^+$, so surely, by 
the spectral theorem,
the inverse of $R_{\psi}$ on $L^2$ is $G(L)$.  It suffices then
to show that $G(L)$, restricted to $L^2 \cap \BB$, has an extension
to a bounded operator on $\BB$.  (Indeed, we would then know that
$m_F(L) G(L) f = G(L) m_F f = f$ for all $f \in L^2 \cap \BB$, so this
would hold for all $f \in \BB$ and $m_F(L)$ would be invertible on $\BB$.)
It suffices then to show that $T = G(L)$ satisfies (i), (ii) and (iii) of
Definition \ref{acceptban}, for some $A_0 > 0$.  Surely (i) is satisfied.\\

First note that $m_F(\lambda)$ is smooth for $\lambda > 0$.  Indeed,
if $V = |H|^2$, then $V \in \mathcal{S}(\RR^+)$, and 
$m_F(\lambda) = \sum_{j \in \ZZ} a^{4j}\lambda^2 V(a^{2j}\lambda)$.  
Since $V$ and all its derivatives are bounded and decay rapidly
at $\infty$, the smoothness of $m_F$ follows at once.

Thus $G \in C^{\infty}((0,\infty))$.  If $\lambda > 0$,
choose $l$ with $a^{2l} \leq \lambda \leq a^{2(l+1)}$.  Surely 
$G(\lambda) = G(a^{-2l}\lambda)$, so for any $k$
\[  |G^{(k)}(\lambda)| = a^{-2kl} |G^{(k)}(a^{-2l}\lambda)| \leq 
a^{2k}\lambda^{-k}M, \]
where $M = \max_{1 \leq \lambda \leq a^2} |G^{(k)}(\lambda)|$.
This shows that $\|\lambda^k G^{k}(\lambda)\|_{\infty} < \infty$
for any $k$.  (ii) and (iii) now follow by the spectral
multiplier theorem of Hulanicki-Stein (\cite{FollandStein82}, Theorem
6.25; see also \cite{Christ}).  (Indeed, by that theorem, $G(L): H^1 
\rightarrow H^1$, so (iii) holds.  Also, in the terminology of
\cite{FollandStein82}, the proof of their Theorem 6.25 shows that
$G(L)$ is given by convolution with a kernel of type $(0,r)$ for any $r$,
so (ii) holds as well.) (a) is therefore established.\\

For (b), note that, since $L^2 \cap \BB$ is dense in $\BB$, and since
$\BB \subseteq \mathcal{S}'$ (continuous inclusion), the operators
$S_{\psi,\psi,b}^{\mathcal{F}}$, acting on $L^2 \cap \BB$,
may be extended to operators on $\BB$,
where they are given by 
\newline
$S_{\psi,\psi,b}^{\mathcal{F}}f = 
\sum_{(j,\gamma) \in \mathcal{F}} <f,\psi_{j,b\gamma}>\psi_{j,b\gamma}$.  
It suffices then to show that for some $C > 0$,
the operator norms on $\BB$ of $S_{\psi,\psi,b}^{\mathcal{F}}$ are all less than $C$,
 for all $\mathcal{F}$.
But, during the proof of Theorem \ref{framebdd} (see the discussion after 
(\ref{eq:50})), we have observed that the operator norms of 
$S_{\psi,\psi,b}^{\mathcal{F}}$ on $L^2$ are uniformly bounded in
$\mathcal{F}$, and that the kernels $K = K_{\psi,\psi,b}^{\mathcal{F}}$ 
satisfy the inequality (\ref{xxxyk}) for some $A_0$ independent of $\mathcal{F}$.
This proves our assertion, by definition of acceptable Banach space. \\

For (c), first take $\BB = L^p$, and say $f \in C_c^{\infty}$.   
By Theorem \ref{framebdd} (a):\\
\ \\
(*) Say $\epsilon_1 > 0$.  There is a finite set 
$\mathcal{F}_1 \subseteq \ZZ \times \Gamma$, such that for any finite set $\mathcal{G}$ with 
$\mathcal{F}_1 \subseteq \mathcal{G} \subseteq \ZZ \times \Gamma$, 
we have
$\|S_{\psi,\psi,b} f - S_{\psi,\psi,b}^{\mathcal{G}}f\|_{\infty}
< \epsilon_1$.\\
\ \\
Moreover, 
since the $K_{\psi,\psi,b}^{\mathcal{F}}$
satisfy (\ref{xxxyk}) uniformly
in $\mathcal{F}$, the argument leading to (\ref{nbfnegQ}) shows that
there is a $C$ such that $|[S_{\psi,\psi,b}^{\mathcal{G}}f](x)| \leq C/|x|^Q$
for all $x$ and all finite $\mathcal{G}$.  These facts imply that 
for any $\epsilon > 0$, and any number $0 < q < Q$, there is a finite set 
$\mathcal{F} \subseteq \ZZ \times \Gamma$, such that for any finite set $\mathcal{G}$ with 
$\mathcal{F} \subseteq \mathcal{G} \subseteq \ZZ \times \Gamma$, 
we have
\begin{equation}\label{lpgq}
|S_{\psi,\psi,b} f - S_{\psi,\psi,b}^{\mathcal{G}}f|< \epsilon g_q. 
\end{equation}
(Here $g_q$ is as in Lemma \ref{gmlem}).  If now also $q$ is also required to satisfy
$q > Q/p$, so that $g_q \in L^p$,

then in (\ref{lpgq}), 
$\|S_{\psi,\psi,b} f - S_{\psi,\psi,b}^{\mathcal{G}}f\|_p < \epsilon\|g_q\|_p$.
The unconditional convergence in $L^p$, for $f \in C_c^{\infty}$,
follows at once.\\
\ \\
Finally, in (c), 
take $\BB = H^1$; then $\mathcal{D} = C_c^{\infty} \cap H^1
= \{f \in C_c^{\infty}: \int f = 0\}$.  
We define a standard 
molecule to be an $L^2$ function $M$ with $\|M\|_2 \leq 1$,
$\int |M(x)|^2 |x|^{Q+1} \leq 1$, and $\int M = 0$.  In 
\cite{CoifmanWeiss2}, it is shown that $M \in H^1$, and further
that there is an $A_0 > 0$ such that $\|M\|_{H^1} \leq A_0$ for all
standard molecules $M$. 

Suppose now $f \in \mathcal{D}$.  Combining (*) above and (\ref{zermom}), we see that
for any $\epsilon > 0$, and any number $0 < q < Q+1$, there is a finite set 
$\mathcal{F} \subseteq \ZZ \times \Gamma$, such that for any finite set $\mathcal{G}$ with 
$\mathcal{F} \subseteq \mathcal{G} \subseteq \ZZ \times \Gamma$, 
we have
\begin{equation}\label{lpgq1}
|S_{\psi,\psi,b} f - S_{\psi,\psi,b}^{\mathcal{G}}f| < \epsilon g_q. 
\end{equation}
If now also $q$ is also required to satisfy $q > Q + \frac{1}{2}$, then
$\max(\|g_q\|_2, [\int |g_q(x)|^2 |x|^{Q+1} dx]^{1/2}) = C_0$, say, is finite.
Thus, in (\ref{lpgq1}), $S_{\psi,\psi,b} f - S_{\psi,\psi,b}^{\mathcal{G}}f$
is $\epsilon C_0$ times a standard molecule, so its $H^1$ norm is less than
$\epsilon C_0 A_0$.  The unconditional convergence in $H^1$, for $f \in
\mathcal{D}$, follows at once.\\

\section{Remarks} 
\begin{enumerate}
\item
When studying frames, 
one often takes several different $\psi$s,
say $\psi^1, \ldots, \psi^N$, all having integral zero, and asks when
$\cup_{k=1}^N \{\psi^k_{j,b\gamma}\}$ is a frame.   
In our situation, by Theorem \ref{genapp} (b),  \[\| \sum_{k=1}^N[\frac{1}{Vb^Q} R_{\psi^k} - S_{\psi^k,\psi^k,b}]\| \leq \frac{C}{b^{Q-1}}. \]

Thus a simple modification of the proof of Theorem \ref{genapp} (c) shows that,
if we can find positive $A,B$ with 
\begin{equation}
\label{eq:62}
AI \leq \sum_{k=1}^N R_{\psi^k} \leq BI, 
\end{equation}
then for some $b_0 > 0$, if $0 < b < b_0$, we can choose $A_b, B_b > 0$ such that 

\begin{equation}\label{eq:63}
A_b\|f\|_2^2  
\leq \sum_{k=1}^N \sum_{j \in \ZZ, \gamma \in \Gamma}|<f,\psi^k_{j,b\gamma}>|^2
\leq B_b\|f\|_2^2 
\end{equation}
for all $f \in L^2$, and such that
\begin{equation}\label{eq:64}
\lim_{b \rightarrow 0^+} \frac{B_b}{A_b} = \frac{B}{A}.
\end{equation}
We restrict attention to $\psi^k$ of the form $\psi^k = F^k(L)\delta
= LH^k(L)\delta$, where $F^k(\lambda) = \lambda H^k(\lambda)$ and 
$H^k, F^k \in \mathcal{S}(\RR^+)$.  Then $\sum_{k=1}^N R_{\psi^k} = 
\sum_{k=1}^N m_{F^k}(\delta)$, and we can take
\begin{equation}\label{eq:62a}
B = \sup_{\lambda > 0} \sum_{k=1}^N m_{F^k}(\lambda),
\ \ \ A = \inf_{\lambda > 0} \sum_{k=1}^N m_{F^k}(\lambda),
\end{equation}
provided this $A$ is positive. 
(This will be so if there does not exist a $\lambda_0
> 0$ such that $F^k(a^{2j}\lambda_0) = 0$ for all $k$ and all integers $j$.)
With higher $N$, one has more flexibility in making $\sum_{k=1}^N m_{F^k}$
nearly constant, thereby getting a nearly tight frame.

\item
In this article, we have let $L$ be the sublaplacian for simplicity, but the
our main results (Theorem \ref{T:2}, Corollary \ref{C:1} and Theorem \ref{T:3})
continue to hold if $L$ is any positive Rockland operator.  (In (\ref{daubcrit}),
one must change $a^{2j}$ to $a^{kj}$ where $k$ is the homogeneous degree of $L$.)
Indeed, the key fact that we have used about $L$ is Theorem \ref{T:4}, and that
theorem continues to hold if $L$ is a positive Rockland operator.  (See 
\cite{Hulanicki84} for this fact and the definition of a Rockland operator.)

\end{enumerate}

\section{Appendix: T(1) theorem technicalities}

As we have said, the ``easier case'' of the $T(1)$ theorem for stratified Lie groups,
as used in this article, may be proved by making only
 minor changes in the proof for $\RR^n$ in \cite{Stein93}, pages 293-300.  
However, one change requires a little thought.

Stein assumes that $T$ is given to be a continuous linear mapping from $\mathcal{S}$
to $\mathcal{S}'$.  He however only uses this assumption in the argument at the
top of page 296.  Moreover the argument at the top of page 296 uses the Fourier transform.
We need to present a replacement for that argument, for general $G$, in
which only $C_c^{\infty}$ functions are used.

For $R > 0$, let $B(0,R) = \{x \in G: |x| < R\}$.  We begin by observing:

\begin{proposition}
Say $T: C_c^{1}(G) \rightarrow L^2(G)$ is linear and restrictedly bounded.  Then $T$ is
continuous from $C_c^{\infty}$ to $L^2$.
\end{proposition}
{\bf Proof}  By definition we need to show that for any compact set $K \subseteq G$,
there exist $C_0, N$ such that $\|Tf\|_2 \leq C_0$ for all $f \in C_c^{\infty}$
with support contained in $K$ and with $\|f\|_{C^N} \leq 1$.  We claim that we can 
always take $N = 1$.  Indeed,  fix $K$ and choose $R > 0$ with $K \subseteq B(0,R)$.
If $f$ is as above, set $F(x) = cf(Rx)$, where $c = \min(1,R^{-a_1},\ldots,R^{-a_n})$.
Then $F$ is a normalized bump function, and $f = (1/c)F^{R,0}$.  Since $T$ is
restrictedly bounded, for somc $C > 0$, $\|Tf\|_2 \leq CR^{Q/2}/c$, as desired.\\

To replace the argument at the top of page 296 in \cite{Stein93}, we now proceed as
follows.  Say $\phi \in C_c^{\infty}(G)$ has support contained in the
unit ball $B(0,1)$.  For $f \in L^2(G)$, let \[ S_j f = f* \phi_{2^{-j}}.  \]
We claim:

\begin{proposition}
Suppose a linear operator $T: C_c^{1}(G) \rightarrow L^2(G)$ 
is restrictedly bounded.  Then:\\
(a) For all $f \in C_c^{\infty}$, $S_jTS_j f \rightarrow Tf$ in $L^2$ as $j \rightarrow \infty$; and\\
(b) For all $f \in C_c^{\infty}$, $S_jTS_j f \rightarrow 0$  in $L^2$ as $j \rightarrow -\infty$.
\end{proposition}
{\bf Proof} Of course $S_j$ is bounded on $L^2$ for all $j$, and $\|S_j\|
\leq \|\phi_{2^{-j}}\|_1 = \|\phi\|_1 = A$, say.\\
For (a), we observe
\begin{eqnarray*}
\|S_jTS_j f - Tf\|_2 &\leq& \|S_jT(S_jf - f)\|_2 + \|S_jTf - Tf\|_2 \\
& \leq & A\|T(S_jf - f)\|_2 + \|S_jTf - Tf\|_2 \rightarrow 0
\end{eqnarray*}
as $j \rightarrow \infty$, since $S_jf \rightarrow f$ in $C_c^{\infty}$, 
$T: C_c^{\infty} \rightarrow L^2$ is continuous, and
$S_jTf \rightarrow Tf$ in $L^2$.\\

For (b) we observe $\|S_j T S_j f\|_2 \leq A\|TS_jf\|_2$, so we need only 
show $TS_j f \rightarrow 0$ in $L^2$.  Write $J = -j$, and note
\[ S_j f = f*\phi_{2^J} = (f_{2^{-J-1}}*\phi_{1/2})_{2^{J+1}} \]
As $J \rightarrow \infty$, $f_{2^{-J-1}}*\phi_{1/2} \rightarrow 
\phi_{1/2}$ in $C_c^{\infty}$, where $c = \int_G f$;
moreover, for $J$ sufficiently large the supports of all these
functions are contained in the unit ball.  Thus, we may choose $C_1$ such that
for $J$ sufficiently large, any one of these functions is $C_1$ times a normalized 
bump function.  But for any function $F$, 
\[ F_{2^{J+1}} = 2^{-J-1} F^{2^{J+1},0}, \]
so $\|TS_jf\|_2 \leq CC_1 2^{-(J+1)/2} \rightarrow 0$ as $J \rightarrow \infty$, as desired.\\
 
Another very small point:  
we have defined a normalized bump function to be a $C^1$ function with support
contained in the unit ball, whose $C^1$ norm is less than or equal to $1$; Stein
assumes in addition that the function is smooth.  But our definition only makes 
the hypotheses of Theorems \ref{tof1} and \ref{tof1q} stronger, so of course the
theorems hold.\\

{\footnotesize{STONY BROOK UNIVERSITY AND TECHNISCHE UNIVERSIT\"AT M\"UNCHEN}}


\begin{thebibliography}{99}
\bibitem{Christ} M. Christ, {\em $L^p$ bounds for spectral multipliers on nilpotent
groups}, Trans. Amer. Math. Soc. 328 (1991), 73-81.
\bibitem{CoifmanWeiss1}  R. Coifman and G. Weiss, {\em Analyse harmonique non-commutative
sur certains espaces homogenes}, Lecture Notes in Math, vol. 242, Springer-Verlag, 
Berlin and New York, 1971.
\bibitem{CoifmanWeiss2}  R. Coifman and G. Weiss, {\em Extensions of Hardy space and
their use in analysis}, Bull. AMS 83 (1977), 569-645.                                             
\bibitem{CorwinGreenleaf89}L. Corwin and F.P. Greenleaf, {\em Representations of Nilpotent
 Lie Groups and Their Applications}, Cambridge University Press, Cambridge,
 1989.
 \bibitem{Daubechies92} I.Daubechies, \textit{Ten Lectures on Wavelets}, SIAM,
Philadelphia, Pennsylvania, 1992.
\bibitem{DGM86} I. Daubechies, A. Grossman, and Y. Meyer, 
\textit{Painless nonorthogonal expansions}, J. Math. Phys. 27 (1986), 1271-1283.
\bibitem{DavidJourne84} G.David and J.L.Journ\'e,
\textit{A boundedness criterion for generalized Calder\'on-Zygmund operators,} 
Ann. Math. 120 (1984) 371-397.
\bibitem{DuffinSchaeffer52} R.J.Duffin, A.C.Schaeffer,
 \textit{A class of nonharmonic Fourier series}, Trans.Amer.Math.Soc. 72 (1952), 341-366.
\bibitem{FeichGro1} H.G. Feichtinger and K. Gr\"ochenig,
\textit{Banach spaces related to integrable group representations and their atomic
decompositions, I}, J. Funct. Anal. 86 (1989), 307-340.
\bibitem{FeichGro2} H.G. Feichtinger and K. Gr\"ochenig,
\textit{Banach spaces related to integrable group representations and their atomic
decompositions, II}, Mh. Math 108 (1989), 129-148.
\bibitem{FollandStein82} G.B.Folland, E.M.Stein, \textit{Hardy spaces on Homogeneous groups},
Mathematical Notes 28, Princeton University Press, 1982.
\bibitem{FJW91} M. Frazier, B. Jawerth, and G. Weiss, \textit{Littlewood-Paley Theory and the 
Study of Function Spaces}, CBMS Reg. Conf. Series in Math., No. 79, Amer. Math. Soc., 
Providence, RI, 1991.
\bibitem{Fuehr02}H.F\"uhr, \textit{Abstract Harmonic Analysis of Continuous 
Wavelet Transforms}, Lecture Notes in Mathematics 1863, Springer Verlag 2005.
\bibitem{Geller80} D.Geller, \textit{Fourier analysis on the Heisenberg group. I. Schwartz
 space}, J. Funct. Anal. 36 (1980), 205--254.
\bibitem{GellerLiou} D. Geller, \textit{Liouville's theorem for homogeneous groups},
Comm PDE 8 (1983), 1665-1677.
\bibitem{GrafakosLi00} L.Grafakos, X.Li, 
\textit{Bilinear operators on homogeneous groups}, J.Operator Theory 
44 (2000), 63-90.
\bibitem{GHLLWW02} J.E. Gilbert, Y.S. Han, J.A. Hogan, J.D. Lakey, D. Weiland, and G. Weiss,
\textit{Smooth Molecular Decompositions of Functions and Singular Integral Operators},
Memoirs of the AMS, Volume 156, \# 742, 2002.
\bibitem{Grochenig91} K. Gr\"ochenig, \textit{Describing functions: atomic decompositions
versus frames}, Mh. Math. 112 (1991), 1-41.
\bibitem{Han00} Y.S. Han, \textit{Discrete Calder\'on-type reproducing formula}, 
Acta Math. Sinica, English Series 16 (2000), 277-294.
\bibitem{Hoermander67} L.H\"ormander, 
\textit{Hypoelliptic second-order differential equations}, Acta Math. 119 (1967), 147-171.
\bibitem{Hulanicki84} A.Hulanicki,\textit{ A functional calculus for Rockland operators on nilpotent 
Lie groups}, Stud. Math. 78, 253--266 (1984).
\bibitem{JerSan} D. Jerison and A. Sanchez-Calle, \textit{Estimates for the heat kernel
for a sum of squares of vector fields}, Indiana Univ. Math. J. 35 (1986), 835-854.
\bibitem{Lemarie84}P.G. Lemari\'e,  \textit{ Alg\`ebres d'op\'erateurs 
et semi-groupes de Poisson sur un espace de nature homog\`ene},
 Universit\'e de Paris-Sud, D\'epartement de Math\'ematiques, Orsay, 1984. iii+130 pp.
\bibitem{Lemarie1} P.G. Lemari\'e, \textit{Base d'ondelettes sur les groupes stratifi\'es},
Bull. Soc. Math. France 117 (1989), 211-232.
\bibitem{Lemarie2} P.G. Lemari\'e, \textit{Wavelets, spline interpolation and Lie groups},
Harmonic Analysis (Sendai, 1990), Springer, Tokyo, 1991, 154-164.  
\bibitem{LiuPeng97} H.Liu, L.Peng, \textit{ Admissible wavelets associated  with the Heisenberg group}, Pacific Journal of
Mathematics, 180 (1997), 101-123.
\bibitem{Maggioni05} M.Maggioni, \textit{ Wavelet frames on groups and hypergroups via 
discretization of Calder\'on formulas}, Monatshefte f\"ur Mathematik, 2005 (in press).
\bibitem{Mayelithesis05} A.Mayeli, \textit{Discrete and continuous wavelet 
transformation on the Heisenberg group}, Ph.D thesis, Technische Universit\"at M\"unchen, 2005. 
\bibitem{Stein93} E.M.Stein, \textit{Harmonic Analysis}, Princeton Mathematical Series 45,
Princeton University Press, 1993.
\bibitem{Varapoulos} N. Varopoulos, \textit{Analysis on Lie Groups}, 
J. Func. Anal. 76 (1988), 346-410.
\end{thebibliography}
\end{document}